\documentclass[12pt]{amsart}

\usepackage{latexsym}
\usepackage[all]{xy}
\usepackage{amsmath, amsthm }
\usepackage{amssymb}
\usepackage{amsfonts}
\usepackage{color}
\usepackage{mathtools}
\usepackage{cancel}
\usepackage[utf8]{inputenc}

\usepackage{hyperref}
\usepackage{xcolor}
\hypersetup{colorlinks=true,urlcolor=blue,linkcolor=blue,citecolor=blue}

\usepackage{cite}

\newcommand {\QCoh} {\mathbf{QCoh}}
\newcommand {\Cris} {\mathbf{Cris}}

\newcommand {\Vect} {\mathbf{Vect}}

\newcommand{\Gm}{\mathbb{G}_{\mathrm{m}}}

\newcommand {\Map} {\mathbf{Map}}
\newcommand {\Parf} {\mathbf{Perf}}

\newcommand {\OO} {\mathcal{O}}

\newcommand {\Fol}{\mathcal{F}ol}

\newcommand {\F} {\mathcal{F}}

\newcommand {\A} {\mathcal{A}}

\newcommand {\T} {\mathbb{T}}

\newcommand{\Q}{\mathcal{Q}}

\newcommand{\ZZ}{\mathbb{Z}}

\newcommand{\LL}{\mathbb{L}}
\newcommand{\VV}{\mathbb{V}}
\newcommand{\WW}{\mathbb{W}}
\newcommand {\Spec} {\mathbf{Spec}}
\newcommand  {\dg}     {\mathbf{dg}}
\newcommand  {\fdg}     {\mathbf{dg}^{\mathrm{fil}}}

\newcommand  {\medg}     {\epsilon\mathbf{dg}^{\mathrm{gr}}}

\newcommand  {\St}   {\mathbf{St}}
\newcommand  {\dSt}   {\mathbf{dSt}}

\newcommand{\s}{\infty}
\newcommand{\HH}{\mathbb{H}}
\newcommand{\cH}{\mathcal{H}}

\newcommand{\gloop}{\mathcal{L}^{gr}}

\theoremstyle{plain}
\newtheorem{thm}{Théorème}[subsection]
\newtheorem{df}[thm]{Définition}
\newtheorem{prop}[thm]{Propositon}
\newtheorem{rmk}[thm]{Remarque}
\newtheorem{cor}[thm]{Corollaire}

\newtheorem{lem}[thm]{Lemme}

\title{Classes caractéristiques des schémas feuilletés}

\author{Bertrand To\"en}
\date{Version préliminaire, Août 2020}

\begin{document}

\maketitle

\begin{abstract} 
Dans ce travail, nous étudions la notion de feuilletages dérivés sur des schémas et schémas dérivés généraux et 
de caractéristiques quelconques. Nous introduisons la filtration de Hodge associée à un feuilletage dérivé, 
qui filtre de manière fonctorielle la cohomologie de de Rham dérivée. 
Nous utilisons cette filtration pour étudier les propriétés
d'annulation des classes caractéristiques de complexes parfaits munis de connexions le long de feuilletages
dérivés. En guise d'application, nous démontrons des extensions à la caractéristique positive, et à valeurs
dans la cohomologie cristalline, du théorème d'annulation
de Bott (voir \cite{bott}) et d'existence de résidus pour les feuilletages avec singularités (voir \cite{babo}).
\end{abstract}

\tableofcontents

\section*{Introduction}

Ce travail s'insère dans la série de travaux sur l'utilisation de techniques dérivées en théorie des feuilletages
algébriques (voir \cite{tv1,tv2}). Nous nous intéressons ici à la notion de classes caractéristiques de feuilletages, sans
pour autant se restreindre au cas de la caractéristique nulle, avec comme objectif d'étendre 
un résultat classique au cas de la caractéristique positive, à savoir le théorème d'existence de résidus
pour les feuilletages avec singularités de \cite{babo}. 

Pour cela, nous introduisons une notion générale de \emph{feuilletages dérivés} au-dessus d'un schéma
$X$ au-dessus d'une base $S$ quelconque. Ces feuilletages sont définis comme certains 
schémas dérivés munis d'une action d'un champ en groupes, noté $\cH$, qui contrôle 
les objets mixtes gradués. De manière plus précise, le groupe $\cH$ est le produit semi-direct
de $\Gm$ et du \emph{cercle gradué} introduit dans \cite{mrt}. Il s'agit d'un objet en groupes
dans la catégories des champs affines de \cite{tosche1}, et dont les représentations sont 
les complexes mixtes gradués: la graduation est induite par l'action de $\Gm$, et 
la structure mixte est elle codée dans l'action du cercle gradué. Ainsi, 
un schéma dérivé muni d'une action de $\cH$ est-il une version non-linéaire
d'un complexe mixte gradué, et permet de coder
un objet de type "complexe de de Rham", muni de ses structures gradués (degrés des formes différentielles), 
mixtes (différentielle de de Rham), et multiplicative (produit des formes différentielles).
Moralement, notre définition \ref{d3} voit ainsi les feuilletages à travers leurs complexes
de de Rham associés (formes différentielles le long des feuilles). 

Nous montrons que notre définition de feuilletages est raisonnable, tout d'abord en montrant qu'elle
redonne, dans le cas lisse, la notion d'algèbroïde de Lie 
(voir la proposition \ref{p1}). Ainsi, 
sur un schéma lisse $X$ au-dessus de $S$, notre définition inclue les 
sous-fibrés $D\subset \Omega_{X/S}^1$ satisfaisant la condition d'intégrabilité usuelle
$d(D) \subset V\wedge \Omega_{X/S}^1$. Nous définissons par ailleurs des images réciproques
de feuilletages dérivés, pour des morphismes $X' \to X$ de $S$-schémas, mais aussi 
pour des changements de bases $S' \to S$. Cette possibilité de jouer avec la base $S$ 
sera cruciale pour notre énoncé d'existence de résidus en caractéristique positive.  

Une des constructions principales de ce travail est la \emph{filtration de Hodge}
associée à un feuilletage dérivé sur un schéma. Il s'agit d'une généralisation du fait bien 
connu qu'un feuilltage sur une variété différentielle définit une filtration naturelle sur le
complexe de Rham de cette dernière. Dans le contexte dérivé général qui nous intéresse, la
construction de cette filtration est plus délicate, et est basée sur la 
déformation vers le fibré normal. Cette filtration de Hodge, associée à un feuilletage
dérivé $\F$ sur un schéma $X$, est alors une filtration multiplicative sur le complexe
de cohomologie de de Rham (dérivé et complété) $C^*_{dR}(X/S)$ de $X$ relativement à $S$. Le degré zéro du 
gradué associé calcule la \emph{cohomologie feuilletée de $X$}, ou encore la cohomologie longitudinale
des formes différentielle le long des feuilles. Les pièces graduées supérieures sont quand à elles
les complexes de cohomologie de de Rham des puissances extérieures du complexe normal au feuilletage $\F$.
En utilisant les propriétés de fonctorialité de la filtration de Hodge on démontre le théorème principal de ce travail, à savoir l'annulation
des classes de Chern en cohomologie feuilletée d'un complexe muni d'une connexion plate le long des
feuilles. Le cas des feuilletages usuels sur des variétés complexes est un énoncé bien 
connu qui se démontre par exemple à l'aide de la théorie de Chern-Weil.

\begin{thm}[Thm. \ref{t1}]\label{ti1}
Soit $\F \in \Fol(X/S)$ un feuilletage dérivé sur un $S$-schéma $X$, et $E$ un complexe parfait muni
d'une connexion plate le long de $\F$. Alors, les classes de Chern $c_i(E)$ de $E$ sont nulles
dans $H^*_{dR}(X/\F)=Gr^0H^*_{dR}(X/S)$, la cohomologie de de Rham feuilletée.
\end{thm}

Ce théorème possède un certains nombres de corollaires. Citons par exemple l'énoncé d'annulation suivant, 
extension du théorème d'annulation de Bott (voir \cite{bott}) au cas de la caractéristique positive. 

\begin{cor}[Cor. \ref{c4}]\label{cti}
Soit $k$ un corps parfait de caractéristique positive et $S=\Spec\, \WW(k)$ avec $\WW(k)$ l'anneau des vecteurs de 
Witt sur $k$. Soit $X$ une variété lisse sur $k$ muni d'un sous-fibré intégrable $D\subset \Omega^1_{X/k}$
de rang $d$. On suppose que le feuilletage défini par $D$ provient, par changement de base, 
d'un feuilletage dérivé sur $X$ au-dessus de $S$ (ce que l'on appelle une "$S$-structure"). Alors, 
pour tout polynôme homogène $\phi \in \WW(k)[X_1,\dots,X_d]$ de degré $q > d$ (avec $deg(X_i)=2i$), 
on a
$$\phi(c_1(D),\dots,c_d(D))=0 \in H^{2q}_{cris}(X/\WW(k)),$$
où $H_{cris}^*(X/\WW(k))$ désigne la cohomologie cristalline de $X$ au-dessus de $\WW(k)$.
\end{cor}

Le corollaire précédent est la version en caractéristique positive du théorème d'annulation de \cite{bott}.
Cependant, le fait que le feuilletage doive possèder une $S$-structure est ici un nouvel ingrédient nécessaire
pour que l'énoncé soit vrai en cohomologie cristalline. Il est important de noter aussi que, bien que $D$ définisse un feuilletage
au sens classique, la notion de $S$-structure utilise de manière cruciale la notion
de feuilletages dérivés, simplement car $X$ n'est pas un $S$-schéma lisse. 
Enfin, le corollaire possède un raffinement pour les feuilletages à singularités, qui 
devient 
alors une version de l'existence des résidus au sens de \cite{babo} en caractéristique positive (Cor.\ref{c5}). \\

\textbf{Remerciements:} Ce travail a été partiellement financé par l'ERC dans le cadre du programme Horizon 2020  
(projet NEDAG ADG-741501). \\

\textbf{Notations et conventions.} On fixe $k$ un anneau commutatif de base et on notera 
$S=\Spec\, k$. Les schémas, champs ou champs dérivés seront pas défaut au-dessus de $S$ (sauf mention
contraire). L'$\s$-catégorie des champs dérivés (au-dessus de $k$) est notée $\dSt_k$. La sous-$\s$-catégorie
des champs tronqués (i.e. non-dérivés) est notée $\St_k \subset \dSt_k$. L'inclusion 
$\St_k \hookrightarrow \dSt_k$ possède un adjoint à droite, la troncation $X \mapsto t_0X$.

Les complexes seront tous $\ZZ$-gradués, avec une différentielle cohomologique (qui augmente le degré de $1$).
Les filtrations seront elles aussi $\ZZ$-graduées, et de plus croissantes. Nous travaillerons régulièrement
avec des complexes gradués, auquel cas, nous réservons le terme de "degré" pour les degrés
cohomologiques, et nous utiliserons "poids" pour désigner la seconde graduation. De manière générale, 
un complexe gradué $E$ sera noté symboliquement $E=\oplus_{n\in \ZZ} E(n)$, où $E(n)$ est pur de poids $n$. 
Enfin, si $E$ est un complexe gradué et $i\in \ZZ$, $E(i)$ désignera le complexe gradué pour lequel
les poids sont décalés de $i$: $(E(i))(n):=E(n-i)$.

Par convention, l'action canonique de $\Gm$ sur la droite affine $\mathbb{A}^1$, sera
de sorte à ce que la graduation induite sur l'anneau des fonctions $k[T]$ coincide
avec la graduation usuelle des polynômes où $T$ est de poids $1$. Ainsi, géométriquement, 
il s'agit de l'action de $\Gm$ sur $\mathbb{A}^1$ de poids $-1$: pour tout corps
$k$, l'action induite de $k^*$ sur $k$ et donnée par $\lambda(x):=\lambda^{-1}x$
pour $\lambda \in k^*$ et $x\in k$. 

\section{Préliminaires et rappels}

\subsection{Champs dérivés équivariants}\label{Gst}

Le contenu de cette sous-section se trouve dans \cite{to1,tosche1,hagII}

Par définition, un \emph{champ dérivé en groupes} est un objet en groupes dans l'$\s$-topos $\dSt_k$. 
Pour un tel groupe $G$ de $\dSt_k$ on dispose de son champ classifiant $BG:=[*/G] \in \dSt_k$, qui est muni
d'un point de base canonique $* \to BG$. La construction $B$ induit un $\s$-foncteur
de l'$\s$-catégorie des champs dérivés en groupes vers celle des champs dérivés pointés. Cet $\s$-foncteur
est pleinement fidèle et son image essentielle consiste en tous les champs
dérivés pointés $F$ tels que le faisceau des composantes connexes soit trivial ($\pi_0(F)\simeq *$). L'$\s$-foncteur
inverse de $B$ est donné par le groupe des lacets $F \mapsto \Omega_*F=*\times_{F}*$. 

Soit $G$ un champ dérivé en groupes et $BG \in \dSt_k$ son classifiant. On dispose de l'$\s$-catégorie
$G-\dSt_k$, des objets $G$-équivariants $\dSt_k$. Tout comme $\dSt_k$, il s'agit d'un $\s$-topos.
On montre de plus qu'il existe une équivalence naturelle d'$\s$-catégories
$$G-\dSt_k\simeq (\dSt_k)/BG$$
entre champs dérivés $G$-équivariants et champs dérivés au-dessus de $BG$. Cette équivalence
envoie un objet $G$-équivariant $X$ sur le champ quotient $[X/G] \to BG$. L'$\s$-foncteur inverse
envoie $X \to BG$ sur $X_*:=X\times_{BG}*$ muni de sa $G$-action naturelle (définie en faisant 
opérer $G$ sur le $G$-torseur universel $*\to BG$). Par la suite, nous identifierons de manière
implicite ces deux $\s$-catégories à travers cette équivalence.

L'$\s$-foncteur d'oubli
$$G-\dSt_k \longrightarrow \dSt_k$$
possède des adjoints à gauche et un droite. L'adjoint à gauche
associe à $X \in \dSt_k$ l'objet $G\times X$, muni de l'action de $G$ par translation 
à gauche sur le premier facteur. L'adjoint 
à droite est quand à lui donné par le hom interne $X \mapsto Map(G,X)$,
où $G$ opère sur lui-même par translation à droite.

Notons enfin l'existence d'un espace classifiant des structures $G$-équivariantes 
sur un champ fixé $X$. Il s'agit, par définition,
de la fibre de l'$\s$-foncteur d'oubli $G-\dSt_k \to \dSt_k$ prise en l'objet $X \in \dSt_k$. Cette fibre
est un $\s$-groupoïde dont l'espace classifiant est par définition l'espace des
structures $G$-équivariantes sur $X$. Il peut aussi se décrire comme le mapping space 
$$Map_{Gp(\dSt_k)}(G,aut(X))\simeq Map_{*/dSt_k}(BG,Baut(X)),$$
dans l'$\s$-catégorie des champs dérivés en groupes, de $G$ vers le champ des autoéquivalences
de l'objet $X$.

\subsection{Graduations et filtrations}

Suivant \cite{sim,mrt,moul}, un champ dérivé gradué est par définition un objet $X \to B\Gm$, 
où de manière équivalente un champ dérivé $\Gm$-équivariant.
Notons de plus $\A:=[\mathbb{A}^1/\Gm]$, le champ quotient de la droite affine par l'action canonique
du groupe multiplicatif\footnote{Notre convention est ici 
que l'action de $\Gm$ sur l'anneau des fonctions $k[T]$ de $\mathbb{A}^1$ 
définisse la graduation usuelle des polynôme où $T$ est de poids $1$. Ainsi, 
géométriquement, il s'agit en réalité de l'action de $\Gm$ sur $\mathbb{A}^1$ 
de poids $-1$ où $\lambda \in \Gm$ opère par multiplication par $\lambda^{-1}$.}. 
D'après \cite{sim,mrt,moul}, un champ dérivé filtré et un objet 
au-dessus de $\A$. Pour un tel objet $X \to \A$ l'objet sous-jacent est la fibre 
$X_1:=1 \times_{\A}X$, où $1=* \to \A$ est le point global induit par l'unité $1 \in k=\mathbb{A}^1(k)$. 
On dit que alors que $X$ définit une filtration sur $X_1$. Le gradué associé est 
par définition $X^{gr} \to B\Gm$, le changement de base
$$X^{gr} := X \times_{\A}B\Gm$$
où $B\Gm \hookrightarrow \A$ est la gerbe résiduelle du champ $\A=[\mathbb{A}^1/\Gm]$ au point $0 \in \mathbb{A}^1(k)$.
L'objet $X^{gr} \to B\Gm$ définit ainsi un champ dérivé gradué, dont l'objet sous-jacent est 
la fibre $X_0=0 \times_{\A}X$ de $X \to \A$ prise en $0=* \to \mathbb{A}^1 \to \A$. De même, on dispose
d'une notion de champs filtrés relatifs à un champ de base $Z$, comme étant les objets
de $\dSt_k/(\A \times Z)$. 

Il est utile de rappeler ici le résultat de \cite{moul}, qui affirme l'existence d'une
équivalence symétrique monoïdale d'$\s$-catégories
$$\QCoh(\A) \simeq \fdg_k,$$
entre les quasi-cohérents sur le champ $\A$ et 
les complexes filtrés de $k$-modules. Ici, $\fdg_k$ est par définition l'$\s$-catégorie
des représentations de l'ensemble ordonné $(\ZZ,\leq)$, des entiers munis de leur ordre naturel.
Ainsi, un objet $E$ de $\fdg_k$ est-il donné par une suite de complexes $F^iE$, avec $i\in \ZZ$,
muni de morphismes $F^iE \to F^{i+1}E$. L'objet sous-jacent à $E$ est par définition 
$F^\s E:=co\lim_i F^iE$, et son gradué associé $\oplus_i F^iE/F^{i-1}E$, où 
$F^iE/F^{i-1}E$ signifie le cône du morphisme $F^{i-1}E \to F^{i}E$. Nous dirons aussi 
que l'objet filtré $E$ est \emph{complet}, si pour tout $i$, le morphisme naturel
dans $\dg_k$
$$F^iE \longrightarrow \lim_{j\leq i}F^iE/F^jE$$
est une équivalence. Les objets complets de $\fdg_k$ forment une sous-$\s$-catégorie pleine
et l'$\s$-foncteur d'inclusion possède un adjoint à gauche évident, à savoir l'$\s$-foncteur
de complétion (qui remplace chaque $F^iE$ par $\lim_{j\leq i}F^iE/F^jE$).

Nous utiliserons souvent le résultat de \cite{moul} dans la situation suivante. Si 
$p : X \to \A$ est un objet filtré, on dispose d'une image direct de complexes quasi-cohérents
$$p_* : \QCoh(X) \longrightarrow \QCoh(\A)\simeq \fdg_k.$$
Ainsi, pour tout quasi-cohérent $E$ sur $X$ on dispose d'un complexe
filtré $p_*(E)$. Moralement, il faut entendre ce complexe $p_*(E)$ comme la cohomologie de $X_1$ à coefficients
dans $E$, muni de la filtration induite par la filtration $X$ sur $X_1$. Ceci n'est cependant 
correct uniquement sous-certaines hypothèses de finitudes qui permettent 
d'affirmer que $p_*$ commute aux changement de bases le long des deux morphismes $0,1 \to \A$. En effet, la compatibilité aux changement de bases de l'$\s$-foncteur $p_*$
est nécessaire afin de pouvoir identifier le complexe sous-jacent au complexe filtré $p_*(E)$ avec
la cohomologie de $X_1$ à coefficients dans $E$ (restreint à l'ouvert $X_1 \subset X$).

\subsection{Déformation vers le fibré normal}

Soit l'immersion fermée $B\Gm \hookrightarrow \A$ donnée par la gerbe résiduelle en $0$. Cette immersion
définit un objet de $\dSt_k/\A$, et donc un champ dérivé filtré. Nous noterons $\Q \in \dSt_k/\A$ cet 
objet, que nous appellerons, non sans ironie, le \emph{point quantique}. Cet objet définit une filtration
sur l'objet vide $\Q_1=\emptyset$, mais dont le gradué associé est 
$$\Q_0\simeq \Spec\, (k\oplus k[1])=0\times_{\mathbb{A}^1}0,$$ 
le spectre de la $k$-algèbre simpliciale libre sur un générateur en degré $-1$ (et de poids $1$). En d'autres termes, 
$\Q$ produit une dégénérescence du vide vers le \emph{point virtuellement vide} $\Spec\, (k\oplus k[1])$, 
qui peut faire penser à la création d'une paire particule/anti-particule par fluctuation quantique. 
De manière plus sérieuse, cette filtration sur le vide produit par suspension la filtration bien connue
sur $S^0=*\coprod *$ dont le gradué associé est le spectre des nombres duaux $k[\epsilon]$. Par double suspension,
elle produit le cercle filtré de \cite{mrt}.

L'objet filtré $\Q$ est en réalité la déformation vers le fibré normal universelle, au sens suivant. Soit 
$f : X \to Y$ un morphisme dans $\dSt_k$, que l'on considère comme un objet de $\dSt_k/Y$. On considère alors
$\Q\times Y$, comme objet filtré de $\dSt_k/Y$, et l'on forme le hom interne au-dessus de $\A \times Y$
$$\Map_{/\A\times Y}(\Q\times Y,\A\times X) \in (\dSt_k)/(\A \times Y).$$
Le champ dérivé $\Map_{\A\times Y}(\Q\times Y,\A\times X)$ définit ainsi un objet filtré au-dessus de $Y$. L'objet 
sous-jacent n'est autre que $\Map_{/Y}(\emptyset,X)=Y$, et $\Map_{\A\times Y}(\Q\times Y,X)$ définit 
donc une filtration sur l'objet $Y$. L'objet sous-jacent au gradué associé
est quand à lui
$$\Map_{/\A\times Y}(\Q\times Y,X)_1 \simeq \Map_{/Y}(\Spec\, (k\oplus k[1])\times Y,X)=: \T_{X/Y}[1],$$
où $\T_{X/Y}[1]$ est le champ tangent décalé de $X$ relativement à $Y$ (et $\Gm$ opère par
homothétie dans les fibres). Lorsque de plus $X$ et $Y$ sont des champs
d'Artin dérivés, on peut alors écrire ce tangent décalé comme
$$\T_{X/Y}[1]=\VV(\mathbb{L}_{X/Y}[-1]),$$
le champ linéaire associé au complexe quasi-cohérent $\mathbb{L}_{X/Y}[-1]$ (voir \cite{hagII}).
Enfin, la projection naturelle $\Q \to \A$ dans $\dSt_k/\A$, induit un morphisme naturel de champs dérivés filtrés 
au-dessus
de $Y$
$$e : \A\times X=\Map_{/\A\times Y}(\A\times Y,\A \times X) 
\longrightarrow \Map_{/\A\times Y}(\Q \times Y,\A \times X).$$
Sur les objets sous-jacents le morphisme $e$ induit le morphisme d'origine $X \to Y$, alors que sur les gradués
associés il s'agit de la section nulle $X \to \T_{X/Y}[1]$ du tangent décalé. Notons ici que la graduation
sur $\T_{X/Y}[1]$, c'est à dire l'action de $\Gm$, n'est autre que l'action par dilatation
de poids $1$ dans les fibres (la section nulle est ici $\Gm$-équivariante).

Avec les constructions précédentes, nous noterons 
$$Def(f):=\Map_{/\A\times Y}(\Q\times Y,\A\times X)$$
que nous considérerons comme un objet de $\dSt_k/\A$, c'est à dire comme un champ 
dérivé filtré. Il est muni d'un morphisme naturel de champs dérivés filtrés
$$e : X \times \A \longrightarrow Def(f).$$

\begin{df}\label{d1}
Le morphisme de champs dérivés filtrés $e : X \times \A \longrightarrow Def(f)$
construit ci-dessus est \emph{la déformation vers le fibré normal 
de $f : X \to Y$}.
\end{df}

Pour terminer, notons qu'il est aussi possible de travailler au-dessus d'un objet de base $Z \in \dSt_k$
quelconque fixé, en utilisant exactement les mêmes constructions. \`A un 
morphisme $X \to Y$ dans $\dSt_k/Z$, nous associons un morphisme d'objets filtrés
au-dessus de $Z$
$$e : X \times \A \longrightarrow Def(f).$$
Sur les objets sous-jacents ce morphisme est le morphisme d'origine $X \to Y$, alors que le gradué associé
est la section nulle $X \to \T_{X/Y}[-1]$. Lorsque $Z=BG$ pour un objet en groupes $G$, ceci nous
permet de considérer une version $G$-équivariante de la déformation vers le fibré normal. 

\section{Feuilletages dérivés sur des schémas}

\subsection{Le cercle gradué et structures graduées mixtes}

Le contenu de cette sous-section est essentiellement tiré de \cite{mrt}. \\

Nous considérons $k<x>$, la $k$-algèbre à puissance divisée libre sur un générateur. Il s'agit d'une
$k$-algèbre de Hopf pour laquelle la comultiplication est induite par la formule usuelle $\Delta(x)=1\otimes x + x\otimes 1$.
Le schéma en groupes correspondant est noté $\HH:=\Spec\, k<x>$. Il s'agit d'un schéma en groupes affine abélien 
sur $S=\Spec\, k$. D'après \cite{mrt,tosche2} $\HH$ s'identifie naturellement au sous-schéma en groupes du 
groupe additif $\WW$ des gros vecteurs de Witt, formé de l'intersection des noyaux de tous les Frobenius. 

La $k$-algèbre de Hopf $k<x>$ est naturellement graduée, par le degré en $x$ et de ce fait 
$\HH$ possède une action naturelle de $\Gm$. Ainsi, le groupe $\Gm$ opère sur le champ
$B\HH$, ce qui définit une graduation sur ce dernier. Le champ gradué $B\HH$ n'est autre que
le gradué associé au cercle filtré $S^{1}_{fil}$ de \cite{mrt}. Le produit semi-direct 
de $\Gm$ par $B\HH$ sera noté $\cH$, de sorte à ce qu'il entre dans une suite exacte, naturellement scindée, 
de champs en groupes
$$\xymatrix{
B\HH \ar[r] & \cH \ar[r] & \Gm.}$$

Par définition, une structure graduée mixte sur un champ dérivé $X\in \dSt_k$ est une
structure $\cH$-équivariante sur $X$ au sens de notre section \ref{Gst}. Ainsi une telle structure
consiste en la donnée de

\begin{enumerate}
\item un objet $[X/\cH] \to B\cH$ de l'$\s$-catégorie comma $(\dSt_k)/B\cH$,

\item une équivalence dans $\dSt_k$
$$[X/\cH] \times_{B\cH}* \simeq X$$
entre la fibre de $[X/\cH] \to B\cH$ au point de base naturel de $B\cH$ et l'objet $X$.

\end{enumerate}

\begin{df}\label{d2}
L'\emph{$\s$-catégorie des champs dérivés gradués mixtes} est 
$$\epsilon-\dSt^{gr}_k:=\cH-\dSt_k\simeq (\dSt_k)/B\cH,$$
l'$\s$-catégorie des champs dérivés $\cH$-équivariants.
\end{df}

Notons que la section naturelle $\Gm \to \cH$ induit, par oubli de structure, un $\s$-foncteur
$$\epsilon-\dSt_k^{gr} \longrightarrow \Gm-\dSt_k,$$
de l'$\s$-catégorie des champs dérivés gradués mixtes vers celle des champs dérivés gradués. Nous serons amenés
à considérer la fibre de cet $\s$-foncteur, prise en un objet $X \in \Gm-\dSt_k$, forme un 
$\s$-groupoïde, qui par définition est l'espace classifiant des \emph{structures mixtes sur $X$ compatibles
avec la graduation}. \\

D'après \cite[Prop. 4.2.3 (ii)]{mrt}, on sait que l'$\s$-catégorie monoïdale symétrique
$\QCoh(B\cH)$ s'identifie à $\medg_k$, l'$\s$-catégories des complexes
gradués mixtes. Nos conventions différent ici de celles de \cite{mrt}, du fait que
la structure mixte diminue le poids de la structure graduée. Pour plus de détails, 
on sait que $B\cH$ est un champ relativement affine au-dessus de $B\Gm$, au 
sens de \cite{tosche1}. Son algèbre de Hopf cosimpliciale de fonctions est 
$k\oplus k[-1]$, l'extension triviale de carré nul de $k$ par $k[-1]$, où l'action de $\Gm$ est 
de sorte à ce que le facteur $k[-1]$ soit de poids $1$. Ainsi, $\QCoh(B\cH)$ s'identifie 
aux comodules gradués sur $k\oplus k[-1]$, ou encore par dualité aux modules
gradués sur $k\oplus k[1]$, où maintenant le facteur $k[1]$ est de poids $-1$ pour l'action
de $\Gm$. Ces modules forment, par définition, l'$\s$-catégorie des 
complexes gradués mixtes, et comme $k[1]$ est de poids $-1$ la structure mixte diminue
bien les poids de $1$. Par la suite, nous utiliserons régulièrement, et 
souvent de manière implicite, l'identification compatible
aux structures monoïdales symétriques
$$\QCoh(B\cH) \simeq \medg_k.$$

\subsection{Feuilletages dérivés et leurs fonctorialités}

Soit $X$ un schéma dérivé (au-dessus de $S=\Spec\, k$). Nous allons définir 
une $\s$-catégorie $\Fol(X/S)$, des feuilletages dérivés sur $X$ relatifs à $S$, 
comme une certaine sous-$\s$-catégorie pleine de champs 
dérivés gradués mixtes au-dessus de $X$ en un certain sens.

Nous commençons par considérer l'espace des lacets gradué de $X$ (relatif à $S$)
$$\gloop(X/S)=\Map(B\HH,X),$$
le champ dérivé des morphismes de $B\HH$ dans $X$. L'objet 
$\gloop(X/S)$ possède une action naturelle de $\cH$, qui par construction opère naturellement 
sur le champ $B\HH$ et donc sur $\Map(B\HH,X)$. De manière équivalente on peut 
considérer l'$\s$-foncteur d'oubli
$$\epsilon-\dSt_k \longrightarrow \Gm-\dSt_k,$$
des champs dérivés gradués mixtes vers les champs dérivés gradués. Cet $\s$-foncteur
possède un adjoint à droite, que l'on peut appliquer à $X$ muni de sa graduation triviale (action 
triviale de $\Gm$), pour obtenir $\gloop(X/S)$.
D'après \cite[Thm. 5.1.3(2)]{mrt}
le champ dérivé $\gloop(X/S)$ est lui-aussi représentable par un schéma dérivé au-dessus de $S$ dont le
tronqué coincide avec celui de $X$. On dispose de plus d'une formule explicite
$$\gloop(X/S)\simeq \Spec_X(Sym_{\OO_X}(\LL_{X/S}[1])=\T_{X/S}[-1],$$
décrivant $\gloop(X/S)$ comme le tangent décalé de $X$ au-dessus de $S$. De plus, la graduation 
naturelle sur $\T_{X/S}[-1]$ est l'action de $\Gm$ par homotéthie (de poids $1$) 
dans les fibres de $\T_{X/S}[-1] \to X$.
Noter que cette action se revêle être l'action de poids $-1$ sur le complexe $\LL_{X/S}[1]$.

Supposons maintenant donnée un champ dérivé gradué mixte $\F$, muni d'un morphisme
d'objets gradués mixtes $\pi : \F \to \gloop(X/S)$. Par adjonction, 
$\pi$ correspond à un morphisme de champs dérivés gradués $p : \F \longrightarrow X$,
où $X$ est muni de la graduation triviale. Supposons de plus que $p$ soit représentable
et affine, de sorte à ce que $\F$ soit aussi un schéma dérivé.
Alors, $\OO_\F:=p_*(\OO)$ est un 
complexe quasi-cohérent $\Gm$-équivariant sur $X$ et en tant que tel se décompose en une somme de composantes
isotypiques
$$\OO_\F \simeq \oplus_{n\in \ZZ}\OO_F(n).$$
Par propriété universelle, l'inclusion naturelle $\OO_F(-1) \longrightarrow p_*(\OO)$
induit un morphisme de schémas dérivés gradués
$$\F \longrightarrow \VV(p_*(\OO_F)(-1))=\Spec_X \,  (Sym_{\OO_X}(p_*(\OO_\F)(-1))).$$

Nous posons alors la définition suivante.

\begin{df}\label{d3}
Avec les notations précédentes, un \emph{feuilletage dérivé au-dessus de $X$ relatif à $S$}
consiste en la donnée d'un champ dérivé gradué mixte $\F$, muni d'un morphisme
d'objets gradués mixtes $\F \to \gloop(X/S)$ et vérifiant les trois
conditions suivantes.
\begin{enumerate}
\item La projection induite $\F \to X$ est représentable et affine.

\item Le complexe quasi-cohérent $\OO_\F(-1)$ est d'amplitude $]-\s,-1]$.

\item Le morphisme naturel de schémas dérivés gradués
$$\F \longrightarrow \VV(\OO_\F(-1))$$
est une équivalence.
\end{enumerate}
\end{df}

Les feuilletages dérivés au-dessus de $X$ forment une $\s$-catégorie $\Fol(X)$, pour laquelle les morphismes
sont les morphismes d'objets gradués mixtes au-dessus de $\gloop(X/S)$: $\Fol(X)$ est ainsi définie
comme la sous-$\s$-catégorie pleine de $\epsilon-\dSt_k^{gr}/\gloop(X/S)$ fromée
des objets vérifiant les conditions $(1)-(3)$ de la définition \ref{d3}.

\begin{df}\label{d3'}
L'$\s$-catégorie des feuilletages dérivés sur $X$ relatifs à $S$ est notée $\Fol(X/S)$.
\end{df}

\begin{rmk}\label{r2}
\emph{La définition précédente de feuilltages dérivés est compatible avec celle donnée dans \cite{tv2}. 
Pour voir cela, il faut noter que l'$\s$-catégorie des schémas affines dérivés
$\cH$-équivariants est, en caractéristique nulle, équivalente à celle
des dg-algèbres commutatives graduées mixtes et connectives sur $k$. Les conditions de la définition
\ref{d3} sont alors équivalentes à celle de \cite[1.2.1]{tv2}.
Notons aussi, qu'il existe une extension possible de la définition \ref{d3} au cas non-connectifs, 
c'est à dire en ne supposant plus que $\F$ soit relativement affine sur $X$ mais seulement 
qu'il soit de la forme $\VV(\OO_\F(-1))$. Nous n'utiliserons pas cette généralisation, 
qui correspond à des feuilletages induits par des actions de $n$-groupoïdes généraux, dans ce texte.}
\end{rmk}

\begin{rmk}\label{r2'}
\emph{L'$\s$-catégorie $\Fol(X/S)$ est fortement sensible à la base $S$, et pour
des morphismes de schémas $X \to S' \to S$ les deux $\s$-catégories
$\Fol(X/S)$ et $\Fol(X/S')$ peuvent être très différentes.}
\end{rmk}

Pour $\F \in \Fol(X/S)$ un feuilletage dérivé (sur un schéma dérivé $X$), considérons
la projection naturelle $p : \F \to X$. Rappelons que cette projection est affine et $\Gm$-équivariante, et 
par ailleurs induit une équivalence de champs dérivés au-dessus de $X$
$$\F \simeq \Spec_X\, (Sym_{\OO_X}(\OO_F(-1))).$$
Le complexe quasi-cohérent $\OO_F(-1)[-1]$ sur $X$ sera appelé le \emph{complexe cotangent de $\F$}. Notons que
la condition $(2)$ de la définition \ref{d3} implique que $\LL_\F$ est 
connectif (i.e. d'amplitude contenue dans $]-\s,0]$).

\begin{df}\label{d4}
Le \emph{complexe cotangent} d'un feuilletage dérivé $\F \in \Fol(X/S)$ est 
défini par 
$$\LL_\F := \OO_F(-1)[-1].$$
Le \emph{complexe conormal} d'un feuilletage dérivé $\F \in \Fol(X/S)$ est 
la fibre du morphisme canonique $\LL_{X/S} \to \LL_{\F}$. Il est noté $\mathcal{N}^*_\F$ de sorte
à ce qu'il existe une suite exacte de fibration de complexes quasi-cohérents sur $X$
$$\xymatrix{
\mathcal{N}^*_\F \ar[r] & \LL_{X/S} \ar[r] & \LL_\F.}$$
\end{df}

La construction $\F \mapsto \LL_\F$ est clairement fonctorielle en $\F$ et fournit un $\s$-foncteur
$$\LL_- : \Fol(X/S)^{op} \longrightarrow \QCoh(X).$$
Cet $\s$-foncteur est conservatif, de sorte à ce qu'un objet $\F$ de $\Fol(X/S)$
puisse être raisonnablement pensé à travers le complexe $\LL_\F$ muni d'une structure additionnelle
de structure graduée mixte sur le schéma dérivé $\VV(\LL_\F[1])$. 

Nous utiliserons les complexes cotangents et conormaux en particulier pour définir des conditions
de lissité variées. La définition ci-dessous regroupe les trois plus importantes.

\begin{df}\label{d5}
Un feuilletage dérivé $\F \in \Fol(X/S)$ sur un schéma dérivé $X$ est \emph{parfait} (resp. \emph{quasi-lisse}, 
resp. \emph{lisse}) si son complexe cotangent $\LL_\F$ est un complexe parfait (resp. un complexe
parfait d'amplitude $[-1,0]$, resp. un fibré vectoriel).

De même, un feuilletage dérivé $\F \in \Fol(X/S)$ sur un schéma dérivé $X$ est \emph{transversalement parfait} 
(resp. \emph{transversalement quasi-lisse}, 
resp. \emph{transversalement lisse}) si son complexe conormal $\mathcal{N}^*_\F$ est un complexe parfait (resp. un complexe
parfait d'amplitude $[-1,0]$, resp. un fibré vectoriel).
\end{df}

Il est important de remarquer que l'$\s$-catégorie $\Fol(X/S)$ possède un objet final, un objet initial et 
des limites arbitraires. L'objet final est clairement $\gloop(X/S)$, muni du morphisme identité 
$\gloop(X/S) \to \gloop(X/S)$. Cet objet sera noté 
$$*_{X/S} \in \Fol(X/S).$$
Notons ici que le complexe cotangent de $*_{X/S}$ n'est autre que 
$\LL_{X/S}$, le complexe cotangent de $X$ relatif à $S$. 
L'objet final de $\Fol(X/S)$ est le morphisme canonique $X \to \gloop(X/S)=\Map(\cH,X)$, 
donné par le \emph{morphisme constant}, adjoint de la projection sur le premier facteur
$X\times \cH \to X$. Son complexe cotangent n'est autre que le complexe nul. Cet objet final 
sera noté
$$0_{X/S} \in \Fol(X/S).$$

Par ailleurs, il est facile de voir que les conditions $(1)-(3)$ de la définition \ref{d3} sont stables
par limites arbitraires, prises dans l'$\s$-catégorie $(\cH-\dSt_k)/\gloop(X/S)$. Ceci montre que
$\Fol(X/S)$ possède des limites arbitraires. Pour un diagramme $\{\F_i\}_{i\in I}$ dans 
$\Fol(X/S)$, le complexe cotangent de l'objet limite $\lim_{i\in I}\F_i$ est donné par la colimite, prise
dans l'$\s$-catégorie comma $\LL_{X/S}/\QCoh(X)$, des complexes cotangents des $\F_i$. En formule
$$\LL_{\lim_{i\in I}\F_i} \simeq co\lim_{i\in I} \LL_{\F_i} \in \LL_{X/S}/\QCoh(X).$$
En particulier, les produits cartésiens existent dans $\Fol(X/S)$, et pour deux feuilletages dérivés
$\F$ et $\F'$ sur $X$, le complex cotangent $\LL_{\F\times \F'}$ entre dans un diagramme
cocartésien de complexes quasi-cohérents
$$\xymatrix{
\LL_{X/S} \ar[r] \ar[d] & \LL_\F \ar[d] \\
\LL_{\F'} \ar[r] & \LL_{\F\times \F'}.
}$$

Notons aussi que les $\s$-catégories $\Fol(X/S)$ sont 
fonctorielles en $X$ par images réciproques. Pour un morphisme de schémas dérivés
$f : X \to Y$ (au-dessus de $S$), on dispose en effet d'un $\s$-foncteur d'image réciproques
$$f^* : \Fol(X/S) \longrightarrow \Fol(Y/S).$$
L'$\s$-foncteur $f^*$ est simplement donné par le changement de base le long 
du morphisme induit $\gloop(X/S) \to \gloop(Y/S)$. Ainsi, si $\F\in \Fol(X)$ est défini par 
un champ dérivé $\cH$-équivariant au-dessus de $\gloop(X/S)$, on a 
$$f^*(\F):=\F\times_{\gloop(X/S)}\gloop(Y/S) \longrightarrow \gloop(Y/S),$$
où le produit fibré est pris dans $\cH-\dSt_k$. Cette formule définit effectivement un $\s$-foncteur $f^* : \Fol(X/S) \to 
\Fol(Y/S)$ car l'on vérifie sans peine que les conditions $(1)-(3)$ de la définition \ref{d3} sont préservées. 
Notons aussi que $\LL_{f^*(\F)}$ entre dans un diagramme cocartésien de complexes quasi-cohérents sur $Y$
$$\xymatrix{
f^*(\LL_{X/S}) \ar[r] \ar[d] & \LL_{Y/S} \ar[d] \\
f^*(\LL_\F) \ar[r] & \LL_{f^*(\F)}.
}$$

Pour terminer, les $\s$-catégories $\Fol(X/S)$ possèdent aussi des fonctorialités en 
le schéma de base $S$, qui étend les images réciproques ci-dessus au cas
d'un morphisme entre paires. Supposons que l'on ait un carré commutatif de schémas dérivés
$$\xymatrix{
X' \ar[r]^-{f} \ar[d] & X \ar[d] \\
S' \ar[r]_-{g} & S,}$$
où $S' \to S$ est un morphisme de schémas affines. On ne supposera pas nécessairement que ce carré
est cartésien. On dispose alors d'un morphisme de champs gradués mixtes (au-dessus de $S$)
$$\gloop(X'/S') \longrightarrow \gloop(X/S).$$
On peut alors, par pull-back le long de ce morphisme, considérer l'$\s$-foncteur de changement de base
$\cH-\dSt_k/ \gloop(X/S) \longrightarrow \cH-\dSt_k/ \gloop(X'/S')$, qui envoie
$\F \to \gloop(X/S)$ sur $\F\times_{\gloop(X/S)}\gloop(X'/S')$. Cet $\s$-foncteur préserve les
conditions $(1)-(3)$ de la définition \ref{d3} et ainsi définit 
un $\s$-foncteur de changement de bases
$$(f,g)^* : \Fol(X/S) \longrightarrow \Fol(X'/S').$$
Note que lorsque $g : S' \to S$ est l'identité, on retrouve les images réciproques 
de feuilletages au-dessus de $S$. Au niveau des complexes cotangents, pour $\F \in \Fol(X/S)$, on dispose
d'un diagramme cocartésien de complexes quasi-cohérents sur $X'$
$$\xymatrix{
f^*(\LL_{X/S}) \ar[r] \ar[d] & \LL_{X'/S'} \ar[d] \\
f^*(\LL_\F) \ar[r] & \LL_{f^*(\F)}.
}$$

Un cas particulièrement intéressant pour nous sera lorsque $X' \to X$ est l'identité. Alors, pour un 
morphisme de $S$-schémas $X \to S'$, on disposera d'un $\s$-foncteur de changement de schémas de bases
$$(id,g)^* : \Fol(X/S) \longrightarrow \Fol(X/S').$$

\subsection{Le cas lisse}

\begin{prop}\label{p1}
Soit $A$ une $k$-algèbre commutative lisse et $M$ un $A$-module projectif et de type fini.
On considère le schéma dérivé affine $X:=\Spec\, (Sym_A(M[1]))=\VV(M[1])$, muni de sa graduation naturelle.
Alors, l'espace classifiant des structures mixtes graduées sur $X$, compatible à la graduation, est 
discret et en bijection avec l'ensemble des différentielles $k$-linéaires et mutiplicatives sur
la $k$-algèbre commutative graduée $\bigoplus_i (\wedge_A^iM)[-i]$.
\end{prop}

\textit{Preuve.} Supposons que $X=\Spec\, (Sym_A(M[1]))$ soit muni d'une action de
$B\HH$ compatible avec sa graduation natuelle (pour laquelle $M[1]$ est de poids $-1$ d'après
nos conventions). On trouve ainsi une action de $\cH$ sur
le complexe $Sym_A(M[1])$, et donc une structure mixte gradué 
sur $\bigoplus_i (\wedge_A^iM)[-i]$ dont la graduation est la graduation naturelle. 
Cette structure mixte est de plus compatible à la 
structure multiplicative car provient d'une action de $B\HH$ sur l'algèbre
simpliciale $Sym_A(M[1])$. Cela fournit un morphisme de l'espace classifiant 
des structures graduées mixtes sur $X$ compatibles à la graduation vers 
l'ensemble des différentielles $k$-linéaires et mutiplicatives sur
la $k$-algèbre commutative graduée $\bigoplus_i (\wedge_A^iM)[-i]$. Il reste à voir que cela est 
une équivalence d'espaces. 

On procède par un argument standard de décomposition de Postnikov comme suit. 
On note $aut^{gr}(X)$ le champ en groupes $\Gm$-équivariant
des autoéquivalences de $X$. L'espace classifiant en question
se trouve être l'espace des morphismes de champs $\Gm$-équivariants
$$Map_{\Gm-\dSt_k}(B^2\HH,Baut^{gr}(X)).$$
Tout d'abord, $B\HH$ étant connexe, tout morphisme de champs en groupes
$B\HH \to aut^{gr}(X)$ se factorise par le sous-champ en groupes $aut_0^{gr}(X)$ formé 
des auto-équivalences qui sont homotopes à l'identité. Comme $A$ est lisse
on voit sans peine que $aut_0^{gr}(X)$  est simplement le groupe des auto-équivalences
induisant l'identité sur les groupes d'homotopie de $Sym_A(M[1])$. Il nous faut donc décrire
l'espace
$$Map_{\Gm-\dSt_k}(B^2\HH,Baut_0^{gr}(X)).$$
Notons aussi que $B^2\HH$ est un champ dérivé tronqué, c'est à dire un objet de $\St_k \subset \dSt_k$, 
et que par conséquent 
le champ $aut^{gr}(X)$ peut-être considéré comme un champ tronqué et que l'espace
des morphismes ci-dessus peut se calculer dans $\Gm-\St_k$, l'$\s$-catégorie
des champs non-dérivés $\Gm$-équivariants sur $k$.

Nous allons procéder par décomposition de Postnikov sur 
$Baut^{gr}_0(X)$. Toud d'abord, $\pi_0(Baut^{gr}_0(X))\simeq \pi_1(Baut^{gr}_0(X),*)\simeq *$. Les faisceaux
d'homotopie $\pi_k:=\pi_k(Baut_0^{gr}(X))$, pour $k>1$ sont des faisceaux abéliens $\Gm$-équivariants qui peuvent 
se décrire explicitement comme suit. Pour $B$ une $k$-algèbre commutative, on note $A_B=A\otimes_k B$ 
et de même $M_B=M\otimes_kB$ (muni de sa structure de $A_B$-module naturelle).
Pour $k$-fixé, et $B$ une $k$-algèbre commutative, 
$\pi_k(B)$ s'identifie aux familles d'applications $k$-linéaires
$$d_i : \wedge_{A_B}^iM_B \to \wedge_{A_B}^{i+k-1} M_B$$
pour $0\leq i$, qui sont compatibles avec la structure multiplicative: c'est à dire que
les $d_i$ définissent une application 
$d : Sym_{A_B}(M_B[-1]) \to Sym_{A_B}(M_B[-1])$ qui est une dérivation d'algèbres commutatives 
graduées de degré $k-1$. Cependant, on n'impose pas
plus de conditions aux $d_i$, par exemple nous n'imposons pas $d^2=0$. La compatibilité 
avec la structure multiplicative implique que $d_0$ et $d_1$ déterminent de manière unique
les $d_i$ pour $i>1$. Un élément de $\pi_k(B)$ est donc donné par une paire
$$d_0 : A_B  \to \wedge^{k-1}M_B \qquad d_1 : M_B \to \wedge^{k}_BM_B,$$
où $d_0$ est une dérivation $B$-linéaire, et $d_1$ est un morphisme
$B$-linéaire qui vérifie la règle de Liebnitz
$$d_1(am)=ad_1(m) +(-1)^{k} d_0(a)\wedge m,$$
pour $a\in A$ et $m\in M$.
L'action de $\Gm$ sur ce faisceau est l'action naturelle
de poids $k-1$ qui agit sur les $d_i$ par $\lambda \mapsto \lambda^{k-1}.d_i$.

Pour être plus précis, le faisceau $\pi_k$ possède une projection naturelle
vers le faisceau quasi-cohérent $B \mapsto B\otimes_k Hom_A(\Omega_{A/k}^1,\wedge^{k-1}M)$, 
donnée par $(d_0,d_1) \mapsto d_0$. Cette projection fait de $\pi_k$ un torseur sous
un second faisceau quasi-cohérent, à savoir $B \mapsto B\otimes_k Hom_A(M,\wedge^{k}M)$. 
Nous noterons cette décomposition par une suite exacte de faisceaux abéliens
$\Gm$-équivariants
$$\xymatrix{
0\ar[r] & Hom_A(M,\wedge^{k}M) \ar[r] & \pi_k \ar[r] & Der(A,\wedge^{k-1}M) \ar[r] & 0,}$$ 
où l'action de $\Gm$ se fait par homothéties de poids $1-k$ (i.e. induite par 
l'action de $\Gm$ de poids $-1$ sur $M$).

Nous sommes maintenant en mesure de procéder à la décomposition de Postnikov. Pour cela, 
on rappelle que la cohomologie de $B\cH$ est naturellement donnée par 
$$H^*(B^2\HH,\OO)\simeq k[u]$$
avec $u$ en degré $2$, et de poids $1$ pour l'action de $\Gm$ (voir \cite[4.1.1]{mrt}). 
Ainsi, pour tout faisceau quasi-cohérent
donné par un $k$-module $N$, sur lequel $\Gm$ opère par son action canonique de poids $1-k$, on trouve
$$H^*(B^2\HH,N)\simeq N[u]\simeq \oplus_i N(i)[-2i],$$
où $(i)$ désigne le twist de Tate pour les modules gradués, qui décale les poids de $i$.
Ainsi, la cohomologie $\Gm$-équivariante est-elle donnée par la seule composante non-nulle $i=k-1$ 
$$H^{2k-2}_{\Gm}(B^2\HH,N)\simeq N \qquad H^{i}_{\Gm}(B^2\HH,N)=0 \;\; si \; i\neq 2k-2.$$
Il est maintenant facile, et laissé au lecteur de procéder à la théorie de l'obstruction le long
de la tour de Postnikov de $Baut^{gr}_0(X)$. On commence par voir que l'espace
$$Map_{\Gm-\St_k}(B^2\HH,(Baut^{gr}_0(X))_{\leq 2})$$
est discret et en bijection avec $H^2(B^2\HH,\pi_2)$, qui n'est autre que l'ensemble
des paires $d_0 : A \to M$ et $d_1 : M \to \wedge^2M$, avec $d_0$ une dérivation et 
$d_1$ vérifiant Liebnitz par rapport à $d_0$. L'obstruction à relever un tel élément à un morphisme
$B^2\HH \to (Baut^{gr}_0(X))_{\leq 3}$ se trouve dans un $H^4_{\Gm}(B^2\HH,\pi_3)$, 
dont les éléments sont des applications $k$-linéaires
$\wedge^iM \to \wedge^{i+2}M$ vérifiant les conditions adéquates. Cette obstruction 
s'annule si et seulement si les $d_i$ vérifient $d^2=0$. Dans ce cas, l'espace des relêvements
est en bijection avec $H^3_{\Gm}(B^2\HH,\pi_3)=0$, et ainsi de suite. 
\hfill $\Box$ \\

La proposition précédente permet de décrire les feuilletages lisses (au sens
de la définition \ref{d5}) au-dessus d'un schéma $X$ lisse sur $S$, en termes
de structures classiques. Pour un tel schéma lisse $X$, on dispose d'une notion 
d'algèbroïde de Lie sur $X$, définie comme suit. Une telle algèbroïde est la donnée d'un couple
$(V,d)$, où $V$ est un fibré vectoriel sur $X$, et $d$ est une différentielle, $k$-linéaire et multiplicative
sur le faisceaux d'algèbres commutatives graduées $\bigoplus_i (\wedge^iV[-i])$. La notion de morphismes
entre de telles algébroïdes est évidente. On a alors le corollaire suivant.

\begin{cor}\label{cp1}
Soit $X$ un schéma lisse sur $S$. Alors, la sous-$\s$-catégorie pleine de $\Fol(X/S)$ formée
des feuilletages dérivés lisses est équivalente à la catégorie des algébroïdes de Lie sur $X$.
\end{cor}

\section{Theorie de de Rham feuilletée}

\subsection{Cohomologie de de Rham feuilletée}

Soit $X$ un schéma dérivé et $\F \to \gloop(X/S)$ un feuilletage dérivé au-dessus de $X$. Par définition
$\F$ est équipé d'une action du champ en groupes $\cH$, et on peut ainsi former le quotient
$[\F/\cH]$. Par définition, la cohomologie de de Rham feuilletée est la cohomologie, à coefficients 
quasi-cohérente, du champ $[\F/\cH]$, ou de manière équivalente la cohomologie $\cH$-équivariante
de $\F$ à coefficients quasi-cohérents. Afin de formaliser ceci, nous commençons par la définition suivante.

\begin{df}\label{d6}
Avec les notations précédentes, soit $f : [\F/\cH] \to S$ la projection canonique. Pour tout 
objet $E \in \QCoh([\F/\cH])$, le \emph{complexe de Rham feuilleté de $X$ le long de $\F$ et à coefficients 
dans $E$} est 
$$C^*_{dR}(X/\F,E):=f_*(E) \in \QCoh(S)\simeq \dg_k.$$
Lorsque $\F=*_{X/S}$ est le feuilletage final et $E=\OO$ on note simplement
$$C^*_{dR}(X/S):=C^*_{dR}(X/\F,E)$$
que l'on appelle le \emph{complexe de cohomologie de de Rham de $X$ relatif à $S$.}
\end{df}

\begin{rmk}\label{r3}
\emph{Le complexe $C^*_{dR}(X/S)$ est souvent appelé} cohomologie de de Rham dérivée de
$X$ sur $S$ \emph{dans la littérature, voir par exemple \cite{bhatt}. Comme nous n'utiliserons pas la
cohomologie de de Rham naïve (i.e. non dérivée) dans cet article, en dehors du cas lisse où 
ces deux notions coïncident, nous ne mentionnons pas l'adjectif} dérivée \emph{dans la
définition précédente. Par ailleurs, la version que nous utiliserons ici est la cohomologie
de de Rham dérivée et} complétée le long de la filtration de Hodge. \emph{Nous renvoyons
à l'appendice \ref{app} pour quelques détails supplémentaires sur les aspects de complétions}.
\end{rmk}

Nous utiliserons de manière essentielle les propriétés de multiplicativité de la cohomologie de de Rham feuilletée. Il s'agit
ici de remarquer que l'$\s$-foncteur
$$f_* : \QCoh([\F,\cH]) \longrightarrow \QCoh(S)$$
est muni d'une structure naturelle lax-monoïdale symétrique, en tant qu'adjoint à droite de l'$\s$-foncteur
monoïdal symétrique $f^*$. Ainsi, $f_*$ s'étend naturellement en un $\s$-foncteur sur les $\s$-catégories
de $E_\s$-algèbres et de modules. En particulier, nous utiliserons que $C^*_{dR}(X/\F,\OO)$
est naturellement muni d'une structure de $E_\s$-algèbre. \\

Les coefficients naturels pour la cohomologie de de Rham feuilletée semblent être tous les quasi-cohérents sur 
$[\F/\cH]$. Cependant, nous allons restreindre cette classe en définissant une sous-$\s$-catégorie pleine
$\Cris(\F) \subset \QCoh([\F/\cH])$ d'objets que nous appellerons des cristaux quasi-cohérents le long de
$\F$. Dans \cite{tv2} cette $\s$-catégorie est notée $\QCoh(\F)$, et pour éviter de la confondre avec l'$\s$-catégorie
des quasi-cohérents sur le schéma dérivé $\F$, nous la noterons $\Cris(\F)$.

Afin de la définir, considérons la projection naturelle $p : \F \to X$. Cette projection n'est 
pas $\cH$-équivariante, mais elle est $\Gm$-équivariante (pour l'action triviale sur $X$). 
On dispose ainsi d'une image réciproque
$$p^* : \QCoh(X) \longrightarrow \QCoh([\F/\Gm]).$$
Cet $\s$-foncteur s'avère être pleinement fidèle, et les objets de son image essentielle 
sont les quasi-cohérents $\Gm$-équivariants sur $\F$ qui sont \emph{libres sur leur composante
de poids $0$ (pour l'action de $\Gm$)}. En clair, si l'on écrit $\F=\VV(\LL_\F[1])$, un objet 
$E \in \QCoh([\F/\Gm])$ se décrit comme un faisceau de $Sym_{\OO_X}(\LL_\F[1])$-modules gradués sur $X$. En tel 
faisceau est dans l'image essentielle de $p^*$ si et seulement si l'inclusion naturelle $E(0) \to E$ de 
sa composante de poids $0$, qui est un morphisme de $\OO_X$-modules, induit une équivalence
$$Sym_{\OO_X}(\LL_\F[1]) \otimes_{\OO_X}E(0) \simeq E.$$
On peut aussi voir cela en utilisant l'adjoint à droite $p_*$ et en remarquant que 
l'unité de l'adjonction $id \Rightarrow p_*p^*$ est toujours une équivalence. 

Rappelons enfin que la section $\Gm \to \cH$ induit un morphisme naturel 
$[\F/\Gm] \to [\F/\cH]$ et par image réciproque un $\s$-foncteur d'oubli
$\QCoh([\F/\cH]) \longrightarrow \QCoh([\F/\Gm])$.

\begin{df}\label{d7}
Avec les notations précédentes, un objet $E \in \QCoh([\F/\cH])$ est un 
\emph{cristal le long de $\F$} si son image dans 
$\QCoh([\F/\Gm])$ appartient à l'image essentiel de $p^* : \QCoh(X) \to \QCoh([\F/\Gm])$
(i.e. est libre sur sa composante de poids $0$).
\end{df}

Les cristaux le long de $\F$ forment une $\s$-catégorie $\Cris(\F)$, qui par définition est une 
sous-$\s$-catégorie pleine de $\QCoh([\F/\cH])$. Elle peut aussi se décrire comme un produit fibré
$$\Cris(\F)\simeq \QCoh([\F/\cH]) \times_{\QCoh([\F/\Gm])}\QCoh(X).$$
La projection naturelle $\Cris(\F) \longrightarrow \QCoh(X)$ fournit la notion d'un quasi-cohérent
sous-jacent à un cristal, et un objet de $\Cris(\F)$ doit être pensé comme à 
un quasi-cohérent sur $X$ muni d'une \emph{connexion plate le long de $\F$}. Nous renvoyons
à \cite{tv2} pour plus de détails sur cette notion dans le cadre de la caractéristique nulle, dans laquelle
nous montrons que les cristaux peuvent aussi se décrire comme des dg-modules sur un faisceau en dg-anneaux
d'opérateurs différentiels adéquat.\\

Les $\s$-catégories $\Cris(\F)$ sont fonctorielles en $\F$ de la manière suivante. Soit 
$$\xymatrix{
X' \ar[d] \ar[r]^-{f} & X \ar[d] \\
S' \ar[r]_-{g} & S
}$$
un carré commutatif de schémas dérivés, avec $S' \to S$ un morphisme de schémas affines. Nous avons défini
une image réciproque $(f,g)^* : \Fol(X/S) \longrightarrow \Fol(X'/S')$. Pour $\F \in \Fol(X/S)$, 
nous avons un morphisme $\cH$-équivariant canonique $\F \times_{\gloop(X/S)}\gloop(X'/S')\to \F$, et donc une image réciproque
quasi-cohérente
$$(f,g)^* : \QCoh([\F/\cH]) \longrightarrow \QCoh([(f,g)^*(\F)/\cH]).$$
Cette image réciproque préserve les sous-$\s$-catégories de cristaux et induit ainsi 
$$(f,g)^* : \Cris(\F) \longrightarrow \Cris((f,g)^*(\F)).$$

\subsection{La filtration de Hodge feuilletée}

Soit $X$ un schéma dérivé et $\F \in \Fol(X/S)$ un feuilletage dérivé sur $X$. Nous allons 
voir que $\F$ permet de définir une filtration sur la $E_\s$-algèbre 
$C^*_{dR}(X/S)$ de cohomologie de de Rham de $X$ relatif à $S$. Pour cela nous
considérons le morphisme $\pi : \F \to \gloop(X/S)$ de schémas dérivés $\cH$-équivariants, et nous
formons sa déformation au fibré normal comme rappelée dans la définition \ref{d1}. Cela fournit 
un schéma dérivé $\cH$-équivariant filtré $Def(\pi)$. On peut en former le quotient
par $\cH$ pour obtenir un champ dérivé filtré
$$p : [Def(\pi)/\cH] \longrightarrow \A.$$
L'image directe $p_*(\OO)$ fournit alors une $E_\s$-algèbre quasi-cohérente sur $\A$, ou en d'autres termes
une $E_\s$-algèbre filtrée. 

\begin{df}\label{d8}
La $E_\s$-algèbre filtrée $p_*(\OO)$ définie ci-dessus sera notée
$$C^*_{dR}(X/S)^{\F}:=p_*(\OO).$$
\end{df}

Il n'est pas difficile de décrire l'objet sous-jacent à $C^*_{dR}(X/S)^{\F}$ ainsi que son 
gradué associé. En effet, si l'on oublie l'action de $\cH$, 
$Def(\pi)$ est un schéma affine dérivé filtré au-dessus de $X$. Le schéma affine
sous-jacent est $\gloop(X/S)=\VV(\LL_{X/S}[1])=\T_{X/S}[-1]$. La filtration est ainsi 
définie par une filtration sur le faisceau de $\OO_X$-algèbres simpliciales commutatives
$Sym_{\OO_X}(\LL_{X/S}[1])$. Cette filtration est elle-même induite en prenant le $Sym$ d'une filtration sur le complexe
$\LL_{X/S}$, ne possèdant qu'un unique cran non trivial, à savoir le morphisme
$$\mathcal{N}_{\F}^* \longrightarrow \LL_{X/S},$$
où $\mathcal{N}_{\F}^*$ est le complexe conormal au feuilletage $\F$ (i.e. la fibre
de $\LL_{X/S} \to \LL_{\F}$, voir la définition \ref{d4}). De manière plus claire, 
on a une filtration $F^*$ définie sur $\LL_{X/S}[1]$ en déclarant que:
\begin{itemize}
\item $F^i(\LL_{X/S}[1])=0$ si $i<-1$ 
\item $F^{-1}(\LL_{X/S}[1]) = \mathcal{N}_\F^*[1]$
\item $F^i(\LL_{X/S}[1]) = \LL_{X/S}[1]$ si $i\geq 0$,
\end{itemize}
avec les morphismes $F^i \to F^{i+1}$ évidents. Comme il s'agit de complexes quasi-cohérents connectifs, nous
pouvons en prendre les $\OO_X$-algèbres simpliciales commutatives libres, et définir ainsi une filtration induite
sur la $k$-algèbre simpliciale commutative $Sym_{\OO_X}(\LL_{X/S}[1])$. Notons que cette filtration est telle que 
$F^i(Sym_{\OO_X}(\LL_{X/S}[1]))=Sym_{\OO_X}(\LL_{X/S}[1])$ dès que $i\geq 0$. Par construction, 
nous savons aussi que son gradué associé est $Sym_{\OO_X}(\LL_{\F}[1]\oplus \mathcal{N}^{*}_{\F}[1])$,
où $\LL_{F}$ est de poids $0$ et $\mathcal{N}_{\F}^*$ est de poids $-1$. Ainsi, la partie
de poids $-i$ n'est autre que $Sym_{\OO_X}(\LL_{\F}[1])\otimes_{\OO_X}(\wedge^i\mathcal{N}^{*}_{\F})[i]$. Comme la
déformation vers le fibré normal est compatible avec l'action de $\cH$, 
$Sym_{\OO_X}(\LL_{\F}[1])\otimes_{\OO_X} \mathcal{N}^{*}_{\F}[i]$ est naturellement muni d'une structure
$\cH$-équivariante en tant que quasi-cohérent $\Gm$-équivariant sur $\F$. En d'autres termes, 
le quasi-cohérent $\wedge^i\mathcal{N}_\F^*$ sur $X$ est naturellement muni d'une structure de cristal
le long de $\F$ au sens de la définition \ref{d7}. \\

En résumé, nous venons de voir les faits suivants.

\begin{prop}\label{p3}
Pour $\F \in \Fol(X/S)$ un feuilletage dérivé sur un schéma dérivé $X$, 
la $E_\s$-algèbre filtrée $C^*_{dR}(X/S)^{\F}$ vérifie les assertions suivantes.
\begin{itemize}
\item Le filtration $F^*$ sur $C^*_{dR}(X/S)^{\F}$ est telle que
$F^i(C^*_{dR}(X/S)^{\F})=F^{i+1}(C^*_{dR}(X/S)^{\F})$ pour tout $i\geq 0$. De plus, on a 
$$F^0(C^*_{dR}(X/S)^{\F})\simeq C^*_{dR}(X/S).$$
En particulier, la $E_\s$-algèbre sous-jacente à $C^*_{dR}(X/S)^{\F}$ est la $E_\s$-algèbre de cohomologie
de de Rham de $X/S$.

\item Les complexes quasi-cohérents $\wedge^i\mathcal{N}^*_\F$ sont munis de structures canoniques
de cristaux le long de $\F$. 

\item Le gradué associé à la $E_\s$-algèbre filtrée $C^*_{dR}(X/S)^{\F}$
est donnée par
$$Gr^*(C^*_{dR}(X/S)^{\F}) \simeq \bigoplus_{i}
C^*_{dR}(X/\F,\wedge^{i}\mathcal{N}^*_\F)[i],$$
où la composante $C^*_{dR}(X/\F,\wedge^{i}\mathcal{N}^*_\F)[i]$ est de poids $-i$ par convention. 

\end{itemize}

\end{prop}

\begin{rmk}\label{r4}
\emph{Pour des feuilletages usuels sur des variétés lisses complexes, la structure 
de cristal sur le fibré conormal $\mathcal{N}^*_\F$ s'appelle parfois} la connexion de Bott.
\end{rmk}

La filtration de Hodge habituelle sur la cohomologie de de Rham s'obtient lorsque 
$\F$ est l'objet final $0_{X/S}$. Cette filtration sera notée 
$$F^i_{Hod}C^*_{dR}(X/S) \to C^*_{dR}(X/S),$$
et appelée \emph{filtration de Hodge absolue} s'il est nécessaire d'éviter certains ambiguités.
Son gradué associé est comme il se doit $\oplus_i C^*(X,\wedge^i\LL_{X/S}[-i])$, 
la cohomologie de Hodge de $X$ relativement à $S$. \\

Enfin, notons que la contruction de la filtration de Hodge associée à un feuilletage
$\F\in \Fol(X/S)$ est fonctorielle. Ainsi, 
si l'on a un diagramme commutatif
$$\xymatrix{
Y \ar[r]^-{f} \ar[d] & X \ar[d] \\
S' \ar[r]_-{g} & S,}$$
et $\F' \to (f,g)^*(\F)$ un morphisme dans $\Fol(Y/S')$, alors 
on dispose d'un morphisme de $E_\s$-algèbres filtrées
$$C_{dR}^*(X/S)^\F \longrightarrow C_{dR}^*(Y/S')^{\F'},$$
dont le morphisme induit sur les $E_\s$-algèbres sous-jacentes n'est autre que 
la fonctorialité usuelle du complexe de de Rham. En d'autres termes, 
la filtration de Hodge associée à un feuilletage est fonctorielle en $\F$
mais aussi en la paire $X\to S$.

\section{Classes caractéristiques}

\subsection{Classes de Chern en cohomologie de de Rham}

On rappelle (voir \cite[2.2.4]{hagII}) l'adjonction
$$i : \St_k \leftrightarrows \dSt_k : t_0,$$
où $t_0$ est l'$\s$-foncteur de troncation et $\St_k$ est l'$\s$-catégorie des champs (non-dérivés). 
L'$\s$-foncteur $i$ est pleinement fidèle et son image essentielle consiste 
en les champs dérivés \emph{tronqués}. Il s'agit précisément des champs dérivés
qui sont extension de Kan à gauche de champs définis sur les schémas affines non-dérivés. En particulier, 
les champs dérivés tronqués sont stables par colimites arbitraires dans $\dSt_k$.

On note $Gl_\s=colim_n Gl_n$, pour les inclusions $Gl_n \hookrightarrow Gl_{n+1}$ standard. Le champ
$\ZZ \times BGl_\s$ est un champ non-dérivé que nous considérons comme un objet de $\St_k$. Par ailleurs, on dispose 
du foncteur de $K$-théorie $K \in \dSt_k$, obtenu par construction de Waldhausen sur le champ 
dérivé des complexes parfaits $\Parf \in \dSt_k$. Il est notable que $K$ et $\ZZ \times BGl_\s$ possède les mêmes morphismes
vers des champs en groupes abéliens, fait bien connu (voir par exemple \cite{gil} où ce fait est utilisé pour définir
les classes de Chern en K-théorie supérieure) mais que nous rappelons ci-dessous.

\begin{prop}\label{p4}
Soit $H$ un champ dérivé en groupes abéliens. Alors, il existe une équivalence naturelle
d'espaces de morphismes
$$Map_{\dSt_k}(K,H) \simeq Map_{\St_k}(\ZZ\times BGl_\s,t_0(H)).$$
\end{prop}

\textit{Preuve.} On commence par remarquer que le morphisme d'adjonction 
$it_0(K) \to K$ est une équivalence, ou en d'autres termes que $K$ est champ dérivé tronqué. 
Pour cela on considère $\Vect$, le champ des fibrés vectoriels, 
que l'on voit comme un champ en monoïdes à l'aide de la somme directe de fibrés. Le champ 
dérivé $\Vect \in \dSt_k$ est tronqué, car s'écrit comme un réunion disjoint de champs classifiants
de groupes lisses $\coprod_n [*/Gl_n]$. L'inclusion canonique $\Vect \longrightarrow \Parf$, induit 
une équivalence après complétion en groupes
$$Gp(\Vect) \simeq K.$$
Or, il est immédiat de vérifier, simplement à l'aide des propriétés universelles, 
que l'$\s$-foncteur de complétion en groupes préserve les objets tronqués. Ainsi, 
$K$ est tronqué car $\Vect$ et donc $Gp(\Vect)$ est tronqué.

Nous utilisons maintenant l'équivalence entre construction $+$ de Quillen 
et construction de Waldhausen de l'espace de $K$-théorie (voir \cite[1.7.1]{wal}). Cet énoncé
affirme l'existence d'une équivalence de champs non-dérivés
$$t_0(K)\simeq Gp(\Vect) \simeq (\ZZ \times BGl_\s)^+,$$
où $+$ fait référence ici à la construction $+$ de Quillen. De plus, par la propriété 
universelle de la construction $+$ on dispose d'une équivalence d'espaces
$$Map_{\St_k}(\ZZ\times BGl_\s,t_0(H)) \simeq Map_{\St_k}((\ZZ\times BGl_\s)^+,t_0(H)).$$
Ainsi, on trouve une équivalence naturelle
$$Map_{\St_k}(\ZZ\times BGl_\s,t_0(H)) \simeq Map_{\St_k}(t_0(K),t_0(H)) \simeq
Map_{\dSt_k}(it_0(K),H)\simeq Map_{\dSt_k}(K,H).$$
\hfill $\Box$ \\

Nous allons appliquer la proposition au cas particulier du champs dérivé
$X \mapsto C^*_{dR}(X/S)$, sous-jacent à 
la cohomologie de de Rham. Pour cela nous rappelons les calculs, bien connus au moins en caractéristique nulle, 
de la cohomologie de de Rham des champs $BGl_n$. 

Soit $\T_{BGl_n}[-1]$, le tangent décalé du champ $BGl_n$, qui est muni d'une structure
naturelle de champs dérivé mixte gradué. Ainsi, la $E_\s$-algèbre de cohomologie de $\T_{BGl_n}[-1]$ à coefficients
dans $\OO$ définie une $E_\s$-algèbre graduée mixte. Cette dernière n'est autre que 
$(Sym_k(gl_n^*))^{hGl_n}.$
Ici $gl_n$ est l'algèbre de Lie de $Gl_n$, et $gl_n^*$ son $k$-module dual, et 
$Sym_k(gl_n^*)$ la $k$-algèbre libre sur $gl_n^*$. Le schéma en groupe $Gl_n$ opère sur
cette dernière par l'action co-adjointe et $(Sym_k(gl_n^*))^{hGl_n}$ désigne la $E_\s$-algèbre de cohomologie
du schéma en groupes $Gl_n$ à coefficients dans icelle. En caractéristique nulle, 
cette cohomologie est triviale en degrés non-nuls, mais dans notre cas $(Sym_k(gl_n^*))^{hGl_n}$ est une
$E_\s$-algèbre cohomologiquement concentrée en degrés $[0,\s[$. Par troncation, on dispose d'un morphisme
canonique de $E_\s$-algèbres
$$H^0((Sym_k(gl_n^*))^{hGl_n})\simeq Sym_k(gl_n^*)^{Gl_n} \longrightarrow (Sym_k(gl_n^*))^{hGl_n}.$$
Ce morphisme est de plus un morphisme de $E_\s$-algèbres graduées mixtes pour la structure mixte triviale
sur le membre de gauche (il s'agit ici d'utiliser la t-structure naturelle sur $\QCoh(B\cH)$). Notons
au passage que, d'après nos conventions, la graduation sur $Sym_k(gl_n^*)$ est de sorte
à ce que $gl_n^*$ soit de poids $-1$. 

Nous pouvons passer aux sections globales, ou encore à la réalisation des objets gradués mixtes au sens de \cite{mrt}
$$|-| : \QCoh(B\cH) \longrightarrow \QCoh(\Spec\, k)=\dg_k.$$
Ceci nous fournit un morphisme de $E_\s$-algèbres (noté CV pour "Chern-Weil")
$$CV :|Sym_k(gl_n^*)^{Gl_n}|\simeq Sym_k(gl_n^*[-2])^{Gl_n} \longrightarrow |(Sym_k(gl_n^*))^{hGl_n}|\simeq C^*_{dR}(BGl_n/S).$$
Ce morphisme est un morphisme de $E_\s$-algèbres filtrées, où la filtration sur le membre de droite est 
la filtration de Hodge, et la filtration sur le membre de gauche n'est autre que 
celle induite par la graduation naturelle. Si $i\leq 0$, sur le $i$-ème étage ce morphisme induit ainsi
$$Sym^{\geq i}_k(gl_n^*[-2])^{Gl_n} \longrightarrow F^{-i}_{Hod}C^*_{dR}(BGl_n/S).$$
Le $i$-ème coefficients de la fonction polynôme caractéristique\footnote{Ici la convention
est de prendre les polynômes caractéristiques inverses: $P_A(T):=Det(Id+TA)$.} définit un élément canonique
pour $i\geq 1$
$$c_i \in Sym^{i}(gl_n^*)^{Gl_n},$$
et par le morphisme de Chern-Weil on dispose ainsi de morphismes canoniques définis dans $\dg_k$
$$c_i : k[-2i] \longrightarrow F^{-i}_{Hod}C^*_{dR}(BGl_n/S) \longrightarrow C^*_{dR}(BGl_n).$$
Ces constructions sont compatibles avec les inclusions $Gl_n \hookrightarrow Gl_{n+1}$, et induisent ainsi 
un morphisme de Chern Weil
$$CV : Sym_k(gl_\s^*[-2])^{Gl_\s} \longrightarrow C^*_{dR}(BGl_\s/S),$$
qui est un morphisme de $E_\s$-algèbres filtrées. Par convention, les sources et but de ce morphisme sont 
obtenus comme limites
$$Sym_k(gl_\s^*[-2])^{Gl_\s} := \lim_n Sym_k(gl_n^*[-2])^{Gl_n}  \qquad 
C^*_{dR}(BGl_\s/S) := \lim_n C^*_{dR}(BGl_n/S).$$
Les classes de Chern $c_i$ deviennent alors des morphismes canoniquement définis dans $\dg_k$
$$c_i : k[-2i] \longrightarrow F^{-i}_{Hod}C^*_{dR}(BGl_\s/S) \longrightarrow C^*_{dR}(BGl_\s/S).$$

Ces morphismes peuvent aussi êtres vus comme des morphismes de champs
$$c_i : BGl_\s \longrightarrow F^{-i}_{Hod}C^*_{dR}(-/S)[2i].$$
Nous appliquons alors la proposition \ref{p4} pour obtenir des morphismes
canoniques de champs dérivés
$$c_i : K \longrightarrow F^{-i}_{Hod}C^*_{dR}(-/S)[2i].$$
En précompposant avec le morphisme canonique $\Parf \longrightarrow K$ nous obtenons enfin la 
définition définitive des classes de Chern des complexes parfaits en cohomologie de de Rham
$$c_i : \Parf \longrightarrow F^{-i}_{Hod}C^*_{dR}(-/S)[2i],$$
comme morphismes dans $\dSt_k$. De manière équivalente nous pouvons considérer
$c_i$ comme un morphisme dans $\dg_k$
$$c_i : k[-2i] \longrightarrow F^{-i}_{Hod}C^*_{dR}(\Parf/S).$$

\begin{df}\label{d9}
Pour $i\geq 1$, le morphisme dans $\dSt_k$
$$c_i : \Parf \longrightarrow F^{-i}_{Hod}C^*_{dR}(-/S)[2i],$$
défini ci-dessus, est la \emph{i-ème classe de Chern}.
\end{df}

Nous pouvons généraliser les classes caractéristique en fixant au préalable un polynôme
$\phi \in k[X_1,\dots,X_i,\dots]$, homogène de degré $q$ avec $deg(X_i)=2i$. Alors, 
$\phi(c_1,\dots,c_i,\dots)$ définit un élément dans $(Sym^{q}_k(gl_\s^*)[-2q])^{Gl_\s}$, et donc, 
à travers le morphsme $CV$ un morphisme de complexes
$$c_\phi : \Parf \longrightarrow F^{-q}_{Hod}C^*_{dR}(-/S)[qi].$$
Il s'agit de la classe caractéristique associé au polynôme $\phi$.

\subsection{Trivialité des classes de Chern de cristaux}

Nous arrivons enfin à notre théorème principal. Soit $X \to S$ un schéma dérivé et $\F\in \Fol(X/S)$ un 
feuilletage dérivé sur $X$. Soit $E \in \Cris(\F)$ un cristal le long de $\F$, tel que 
le quasi-cohérent sous-jacent soit un complexe parfait sur $X$. Notons 
$E(0)$ ce complexe parfait, qui, on le rappelle, est la partie de poids $0$ de $E$
vu comme quasi-cohérent $\Gm$-équivariant sur $\gloop(X/S)$. 

\begin{thm}\label{t1}
Avec les notations ci-dessus, pour tout $i>0$, le morphisme
$$c_i(E(0)) : k[-2i] \longrightarrow C^*_{dR}(X/S)$$
possède une factorisation canonique
$$k[-2i] \longrightarrow F^{-1}C^*_{dR}(X/S)^\F \longrightarrow C^*_{dR}(X/S).$$
\end{thm}

\textit{Preuve.} Soit $E(0) : X \longrightarrow \Parf$ le morphisme de champs dérivés correspondant 
au complexe parfait $E(0)$. On considère la projection naturelle $p : \F \to X$.
Par définition, la structure de cristal le long de $\F$ sur $E(0)$ est donnée par 
une structure $\cH$-équivariante sur le morphisme composé
$$\xymatrix{\F \ar[r]^-p & X \ar[r]^-{E(0)} & \Parf,}$$
où $\cH$ opère sur $\Parf$ de manière triviale. En d'autres termes, le morphisme $E(0)$ se promeut en 
un morphisme de champs dérivés feuilletés
$$(X,\F) \longrightarrow (\Parf,0_\Parf),$$
où $0_{\Parf}\in \Fol(\Parf)$ est l'objet initial. Par fonctorialité de la filtration de Hodge, on en déduit
par image réciproque
un moprhisme de $E_\s$-algèbres filtrés
$C^*_{dR}(\Parf/S) \longrightarrow C^*_{dR}(X/S)^\F,$
où le membre de droite est filtré par la filtration de Hodge usuelle (i.e. celle
associée au feuilletage initial). Ceci implique, que l'image réciproque le long de $E(0)$ 
$C^*_{dR}(\Parf/S) \longrightarrow C^*_{dR}(X/S)^\F$
vient équipée avec une factorization
$F^{-1}_{Hod}C^*_{dR}(\Parf/S) \longrightarrow F^{-1}C^*_{dR}(X/S)^\F.$
Comme les classes de Chern $c_i$ sont des morphismes $k[-2i] \to F^{-1}_{Hod}C^*_{dR}(\Parf/S)$, 
le théorème s'en suit.
\hfill $\Box$ \\

\begin{rmk}\label{rt1}
\emph{Le théorème \ref{t1} affirme plus que l'existence d'une factorisation, sa preuve montre que la
donnée d'une structure de cristal de long de $\F$ sur le complexe parfait $E(0)$ définit une
factorisation canonique. On pourrait aller plus loin et montrer que l'on construit ici un $\s$-foncteur
de l'$\s$-groupoïde des structures de cristaux sur $E(0)$ vers l'$\s$-groupoïde
des factorisations comme dans la conclusion de \ref{t1}.}
\end{rmk}

En termes de classes de cohomologie, le théorème \ref{t1} possède l'implication suivante: l'annulation
des classes de Chern de cristaux en cohomologie de de Rham feuilletée.

\begin{cor}\label{ct1}
Avec les mêmes notations et conditions que le théorème \ref{t1}, toutes les images des classe
de Chern $c_i(E(0)) \in H^{2i}_{dR}(X/S)$ par la projection naturelle $H^{2i}_{dR}(X/S)\to H^{2i}_{dR}(X/\F)$
sont nulles.
\end{cor}

Notons cependant que le théorème \ref{t1} est plus fin que le corollaire précédent. En effet, 
il affirme que les classes de Chern $c_i(E(0)) : k[-2i] \to C^*_{dR}(X/\F)$ viennent 
équipée d'une homotopie canonique au morphisme nul. Une des conséquence de ce fait va être l'existence
de résidus pour les feuilletages dérivés, extension des résidus de feuilletages singuliers complexes
de \cite{babo}. Pour commencer nous avons la conséquence suivante du dernier corollaire, 
qui est une extension du théorème d'annulation de Bott (voir \cite{bott}).

\begin{cor}\label{c2}
Soit $\F \in \Fol(X/S)$ un feuilletage dérivé transversalement lisse (voir \ref{d5}). On suppose
que son fibré conormal $\mathcal{N}_\F^*$ est de rang $d$ (nous dirons alors que $\F$ est de plus
de codimension $d$). Soit $\phi \in k[X_1,\dots,X_i,\dots]$ un polynôme avec
$X_i$ de degré $2i$, et $\phi$ homogène de degré $q >d$. Alors, pour tout crystal $E \in Cris(\F)$, 
avec $E(0)$ parfait sur $X$, on a
$$\phi(c_1(E),\dots,c_i(E),\dots)=0 \in H^{2q}_{dR}(X/S).$$
\end{cor}

\textit{Preuve du corollaire \ref{c2}.} Il suifft de remarquer que pour
$\F\in \Fol(X/S)$ transversallement lisse et de codimension $d$, la filtration
de Hodge associée à $\F$ vérifie
$$F^{-i}C^*_{dR}(X/S)^\F=0 \quad si \; i>d.$$
En effet, cet énoncé est local sur $X$, et on peut donc supposer que $X=\Spec\, A$
est affine. Au niveau des complexes cotangents on a suite exacte de fibration
$$\xymatrix{
\mathcal{N}^*_\F \ar[r] & \LL_{X/S} \ar[r] & \LL_{\F}.
}$$
On peut représenter cette suite par une suite exacte de cofibration de $A$-modules simpliciaux
cofibrants
$$
\mathcal{N}^*_\F \hookrightarrow \LL_{X/S} \twoheadrightarrow \LL_{\F},$$
et supposer que $\mathcal{N}^*_\F$ est un $A$-module libre de rang $d$.
Dans ce cas, la filtration induite sur $Sym_{A}(\LL_{X/S}[1])$ est telle que
$F^{-i}Sym_{A}(\LL_{X/S}[1])$ soit l'idéal simplicial engendré par 
$Sym^{\geq i}(\mathcal{N}^*_{\F}[1])\simeq \oplus_{j\geq i} \wedge^j\mathcal{N}_\F^*[j]$, 
qui est donc nul pour $i>d$.
\hfill $\Box$ \\

Pour l'existence de résidus, fixons $X$ un schéma dérivé et $E$ un complexe parfait
sur $X$. On se fixe un ouvert $U\subset X$, un feuilletage dérivé $\F_U \in \Fol(U/S)$
transversalement lisse de codimension $d$,
et une structure de cristal sur $E_U$ le long de $\F_U$.
On note $Z=X-U \subset X$ muni de sa structure de sous-schéma
réduit.  On note enfin $H^*_{dR,Z}(X/S)$ la cohomologie 
de de Rham de $X$ à supports dans $Z$, définie comme la cohomologie de la fibre 
du morphisme de restriction
$C^*_{dR}(X/S) \to C^*_{dR}(U/S)$. 

\begin{cor}\label{c3}
Avec les notations et conditions précédentes, soit $\phi \in k[X_1,\dots,X_i,\dots]$ un polynôme avec
$X_i$ de degré $2i$, et $\phi$ homogène de degré $q >d$. Alors,
il existe un élément canoniquement défini
$$Res_{\phi}(E) \in H^{2q}_{dR,Z}(X/S)$$
tel que son image dans $H^{2q}_{dR}(X/S)$
soit égale à $\phi(c_1(E),\dots,c_i(E),\dots)$.
\end{cor}

\textit{Preuve du corollaire \ref{c3}.} Notons que le corollaire \ref{c2} implique que
des éléments $Res_\phi$ avec la propriété cherchée existent bien. Cependant, il est question
d'en construire un plus beau que les autres. Cet élément cherché est obtenu à l'aide du théorème
\ref{t1} appliqué à $U$ et $E_U$. En effet, le théorème implique que le morphisme induit par $\phi$ à l'aide
du morphisme de Chern-Weil (voir après la définition \ref{d9})
$$c_\phi(E_U):=\phi(c_1(E_U),\dots,c_i(E_U),\dots) : k[-2q] \longrightarrow C^*_{dR}(U/S)$$
se factorise naturellement par $F^{-q}C^*_{dR}(X/S)^{\F}$. Or, comme le feuilletage
$\F_U$ est transversallement lisse de codimension $d<q$, nous avons déjà vu lors de la preuve
du corollaire \ref{c2} que $F^{-i}C^*_{dR}(X/S)^{\F}$ pour tout $i>d$. Ainsi, $c_\phi$ vient avec
une factorisation canonique au morphisme nul dans l'$\s$-catégorie $\dg_k$. Cette factorisation
induit un morphisme $Res_\phi$ bien défini, de sorte à avoir un diagramme commutatif
dans $\dg_k$
$$\xymatrix{
C^*_{dR,Z}(X/S) \ar[r] & C^*_{dR}(X/S) \ar[r] & C^*_{dR}(U/S) \\
 & k[-2q] \ar[lu]^-{Res_\phi} \ar[u]_-{c_\phi(E)} \ar[ru]_-{c_\phi(E_U)}}$$
\hfill $\Box$ \\

Le corollaire \ref{c3} s'applique en particulier au cas où $\F$ est défini
sur $X$ tout entier, et où $E=\mathcal{N}^*_{\F}$. On trouve alors l'existence
de résidus associés aux feuilletages dérivés qui sont génériquement transversallement 
lisses, dont l'énoncé est plus proche du théorème principal de \cite{babo}.

\begin{cor}\label{c3'}
Soit $\F\in \Fol(X/S)$ un feuilletage dérivé parfait et qui est transversallement lisse de 
codimension $d$ sur un ouvert $U\subset X$. Alors, pour tout $\phi \in k[X_1,\dots,X_i,\dots]$ un polynôme avec
$X_i$ de degré $2i$, et $\phi$ homogène de degré $q >d$, il existe un élément canonique
$$Res_{\phi}(\F) \in H^{2q}_{dR,Z}(X/S)$$
tel que son image dans $H^{2q}_{dR}(X/S)$
soit égale à $\phi(c_1(\mathcal{N}^*_\F),\dots,c_i(\mathcal{N}^*_\F),\dots)$.
\end{cor}

\subsection{Une formule des résidus à la Baum-Bott en caractéristique positive}

Pour terminer, nous allons appliquer les résultats précédents au cas particulier 
où $X$ est une variété lisse sur un corps parfait $k$ de caractéristique $p>2$.
Nous noterons $S=\Spec\, \WW(k)$, l'anneau des vecteur de Witt $p$-typiques sur $k$, et 
$i : s=\Spec\, k \hookrightarrow S$ son unique point fermé. Nous introduisons alors la notion suivante de $S$-structures
sur des feuilletages sur $X$ relatifs à $s$, qui décrit le fait que les "feuilletages sur $X/s$ peuvent 
s'étendre dans la direction de $S$". 

\begin{df}\label{d10}
Soit $\F\in \Fol(X/s)$ un feuilletage dérivé sur la variété $X$ relativement à $s$. Une
\emph{$S$-structure sur $\F$} est la donnée de $\F' \in \Fol(X/S)$ et d'une équivalence
$i^*(\F')\simeq \F$ dans $\Fol(X/s)$.
\end{df}

Un des intérêt de la notion de $S$-structures sur un feuilletage $\F$, et de pouvoir
étendre la filtration de Hodge de $\F$, définie sur $C^*_{dR}(X/s)$, en une filtration 
sur $C^*_{dR}(X/S)$. Par exemple, supposons qu'il existe un carré cartésien de schémas dérivés
$$\xymatrix{
X \ar[d] \ar[r] & Y \ar[d] \\
s \ar[r] & S,
}$$
de sorte à ce que $\F$ soit obtenu par image réciproque d'un feuilletage $\F_Y\in \Fol(Y/S)$. Alors, 
$\F_Y$ définit une $S$-structure sur $\F$, en choisissant $\F'$ comme étant l'image 
de $\F_Y$ par le $\s$-foncteur image réciproque $\Fol(Y/S) \to \Fol(X/S)$. Cependant, 
la notion de $S$-structure garde un sens sans avoir besoin de l'existence de $Y$. En d'autres termes, 
un feuilletage $\F \in \Fol(X/s)$ peut s'étendre au-dessus de $S$ sans que $X$ s'étende au-dessus
de $S$ pour autant.

Tous les feuilletages dérivés ne possèdent pas de $S$-structures, posséder une $S$-structure est ainsi
une condition non-triviale. \`A ce sujet, les deux 
exemples des objets finaux et initiaux de $\Fol(X/S)$ sont instructifs. Pour commencer
$*_{X/s}$ possède toujours une $S$-structure, à savoir $*_{X/S}$. Ceci se voit simplement en remarquant
que les images réciproques préservent les limites et donc l'objet final. En revanche, 
$0_{X/s}$ ne possède pas toujours une $S$-structure. En effet, supposons qu'une telle $S$-structure
$\F' \in \Fol(X/S)$ existe, et pour simplifier supposons $X=\Spec\, A$ affine et non-dérivé. Alors
$C^*_{dR}(X/\F')$ serait alors une $E_\s$-algèbre $A'$, $\WW(k)$-linéaire, et muni d'un isomorphisme
de $k$-algèbres
$$A'\otimes_{\WW(k)}k \simeq C^*_{dR}(X/s)\simeq A.$$
Ceci impliquerait automatiquement que $A'$ est une $\WW(k)$-algèbre plate, et ainsi que $X$ possède un
relèvement à $S$, ce qui n'est pas vrai en général. De manière générale, 
l'existence d'une $S$-structure sur $0_{X/s}$ est essentiellement équivalente à l'existence d'un 
relèvement de $X$ à un $S$-schéma. Il est bon de garder l'image intuitive que la donnée d'une
$S$-structure sur $\F \in \Fol(X/s)$ consiste essentiellement en la donnée d'un relèvement 
à $S$ de "l'espace des feuilles de $\F$", bien que ce dernier n'existe pas dans la catégorie algébrique. 

L'intérêt, pour nous, des $S$-structures est basé sur la comparaison entre
cohomologie de de Rham  (dérivée) et cohomologie cristalline, dans l'esprit des
résultats de \cite{bhatt}. Nous n'aurons pas besoin de résultats de comparaisons, 
mais simplement de l'existence d'un morphisme naturel de $E_\s$-algèbres $\WW(k)$-linéaires
$$\psi : C^*_{dR}(X/S) \longrightarrow C^*_{cris}(X/S),$$
où le membre de droite est défini comme étant
$$C^*_{cris}(X/S):=\lim_n \mathbb{R}\Gamma((Cris(X/S_n),\OO_{X/S_n}),$$
avec $S_n=Spec\, \WW_n(k)$, et $Cris(X/S_n)$ le topos cristallin de $X$ au-dessus de $S_n$
muni de son faisceau structural $\OO_{X/S_n}$ usuel.
L'existence de ce morphisme est certainement bien connu. Cependant, 
\cite{bhatt} construit ce morphisme uniquement pour la cohomologie de de Rham
dérivée \emph{non-complétée}, alors que par définition $C^*_{dR}(X/S)$
est la cohomologie de de Rham dérivée et complétée pour la filtration de Hodge. 
Fort heureusement, les phénomènes de complétions disparaissent dans notre situation
car $X$ est supposée lisse sur $s$. Les détails se trouvent en appendice.

Nous disposons donc d'un morphisme naturel
$$\psi : C^*_{dR}(X/S) \longrightarrow C^*_{cris}(X/S).$$
Ce morphisme est fonctoriel en $X$, et de plus compatible aux classes de Chern. Ainsi, pour tout
complexe parfait $E$ sur $X$, 
$\psi(c_i(E))$ est égal à la i-ème classe de Chern de $E$ en cohomologie cristalline. 
Nos théorème d'annulation en cohomologie de de Rham impliquent donc les énoncés analogues en 
cohomologie cristalline. 
Le résultat suivant est un exemple d'application
du corollaire \ref{c3'} de la section précédente, que nous proposons comme version du théorème
d'annulation de Bott de \cite{bott} au-dessus de corps de caractéristique $p>0$.

\begin{cor}\label{c4}
Soit $X$ une variété lisse sur $k$, et $D\subset \Omega_{X/s}^1$ un sous-fibré de rang $d$ qui 
vérifie la condition d'intégrabilité $d(D) \subset D\wedge \Omega_{X/s}^1$. 
On suppose que l'objet $\F_D \in \Fol(X)$ définit par $D$ (voir \ref{cp1}) possède une $S$-structure 
au sens de \ref{d10}. Alors, pour tout polynôme
$\phi \in \WW(k)[X_1,\dots,X_d]$  avec
$X_i$ de degré $2i$, et $\phi$ homogène de degré $q >d$, on a 
$$\phi(c_1(D),\dots,c_d(D))=0 \in H^{2q}_{cris}(X/\WW(k)).$$
\end{cor}

Les corollaires \ref{ct1}, \ref{c2} et \ref{c3'} possèdent eux-aussi des applications évidentes. Par exemple, 
nous déduisons l'énoncé suivant d'existence de résidus pour des feuilletages singuliers 
en caractéristique positive.

\begin{cor}\label{c5}
Soit $X$ une variété lisse sur $k$, et $D\subset \Omega_{X/s}^1$ un sous-faisceau cohérent 
vérifiant la condition d'intégrabilité $d(D) \subset D\wedge \Omega_{X/s}^1$. On suppose, 
que $D$ est un sous-fibré de rang $d$ sur un ouvert $U \subset X$, et on suppose donnée une $S$-structure
sur le feuilletage que $D$ définit sur $U$ relatif à $s$. On note 
$Z=X-U$ muni de sa structure réduite.
Alors, pour tout $\phi \in \WW(k)[X_1,\dots,X_i,\dots]$ un polynôme avec
$X_i$ de degré $2i$, et $\phi$ homogène de degré $q >d$, il existe un élément canonique
$$Res_{\phi} \in H^{2q}_{cris,Z}(X/\WW(k))$$
tel que son image dans $H^{2q}_{cris}(X/\WW(k))$
soit égale à $\phi(c_1(E),\dots,c_i(E),\dots)$.
\end{cor}

Notons que dans le corollaire précédent, le résidu $Res_\phi$ dépend à priori 
du choix de la $S$-structure sur le feuilletage sur $U$ défini par $D$.

\begin{appendix}

\section{Cohomologie de de Rham dérivée et cohomologie cristalline}\label{app}

Soit $k$ un corps parfait de caractéristique $p>2$, $\WW(k)$ son anneau des vecteurs de Witt
et $S=\Spec\, \WW(k)$. Le but de cet appendice est de
rappeler l'existence d'un morphisme $\psi : C^*_{dR}(X/S) \to C^*_{cris}(X/S)$, 
depuis la cohomologie de de Rham d'une variété lisse $X$ sur $k$, vers sa cohomologie
cristalline, fonctoriellement en $X$. Le contenu de cette partie est principalement tirée
de \cite{bhatt}. \\

Tout d'abord, $C^*_{cris}(X/S)$ est défini par la limite des $C^*_{cris}(X/S_n)$, où 
$S_n=\Spec\, \WW_n(k)$
$$C^*_{cris}(X/S):=\lim_n C^*_{cris}(X/S_n).$$
Nous allons définir le morphisme $\psi$ comme le composé
$$C^*_{dR}(X/S) \to \lim_nC^*_{dR}(X/S_n) \to \lim_n C^*_{cris}(X/S_n)=C^*_{cris}(X/S).$$
Ainsi, il nous suffit de construire les morphismes
$$\psi_n : C^*_{dR}(X/S_n) \to C^*_{cris}(X/S_n)$$
de manière compatible en $n$. De plus, si l'on dispose des morphismes $\psi_n$ lorsque
$X$ est affine, on dispose par recollement des morphismes $\psi_n$ pour toute variété $X$. 
Nous allons donc nous concentrer sur le cas des variétés affines lisses au-dessus $k$, 
et des cohomologies de de Rham et cristallines au-dessus de $S_n$ pour $n$ fixé.
La globalisation est laissée au lecteur.

Soit $\mathcal{C}_n$ l'$\s$-catégorie dont les objets sont les $\WW(k)_n$-algèbres simpliciales
commutatives qui sont des rétractes de complexes cellulaires finis (voir \cite[\S 2.1]{tova}). En particulier, 
les objets $A$ de $\mathcal{C}_n$ sont tels que chaque $\WW_n(k)$-algèbre de m-simplexes $A_m$
est une $\WW_n(k)$-algèbre lisses sur $k$. \`A $A\in \mathcal{C}_n$ nous
associons un objet simplicial dans $\fdg_{\WW_n(k)}$, c'est à dire un complexe filtré 
simplicial de $\WW_n(k)$-module, à savoir le complexe de de Rham 
$$dR(A):=[m] \mapsto \Omega_{A_m/\WW_n(k)}^*[-*],$$
où chaque $\Omega_{A_m/\WW_n(k)}^*[-*]$ est muni de sa différentielle
usuelle et de sa filtration de de Rham. Par normalisation on obtient un complexe 
filtré
$N(dR(A)) \in \fdg_{\WW_n(k)}.$
La construction $A \mapsto N(dR(A))$ définit un $\s$-foncteur
$$N(dR(-)) : \mathcal{C}_n \longrightarrow \fdg_{\WW_n(k)}.$$
Tout objet $A \in \mathcal{C}_n$ possède une augmentation
$A \to \pi_0(A)$, où $\pi_0(A)$ est vu comme un objet simplicial constant. Le noyau de cette 
augmentation définit un idéal simplicial $I \subset A$. Soit $D(I)$ l'enveloppe à puissance
divisé de $I$ dans $A$, de sorte à ce qu'il existe un diagramme
commutatif d'anneaux simpliciaux
$$\xymatrix{A \ar[rd] \ar[r] & D(I) \ar[d] \\
 & \pi_0(A).}$$
On sait que $D(I)$, en tant que $A$-algèbre, est muni d'une connexion plate dont on 
peut former le complexe de de Rham
$$dR(D(I)) : [m] \mapsto D(I_m)\otimes_{A_m}\Omega^*_{A_m/\WW_n(k)}[-*].$$
L'objet $dR(D(I))$ est encore un objet simplicial dans $\fdg_{\WW_n(k)}$: pour $m$ fixé
la différentielle de $D(I_m)\otimes_{A_m}\Omega^*_{A_m/\WW_n(k)}$ est encore la différentielle de de Rham, 
et la filtration est la filtration de Hodge $PD$-adique de \cite[V-2.3.1]{bert}, convolution
de la filtration PD-adique sur $D(I_m)$ et la filtration de Hodge. On rappelle que
le $(-i)$-ème cran de cette filtration est donnée par le sous-complexe
$$I_m^{[i]} \to I_m^{[i-1]}\otimes_{A_m}\Omega^1_{A_m/\WW_n(k)} \to 
I_m^{[i-2]}\otimes_{A_m}\Omega^2_{A_m/\WW_n(k)} \to\dots,$$
où $I_m^{[k]}$ désigne la $k$-ème  PD-puissance de $I_m$ dans $D(I_m)$
(avec $I_m^{[k]}=D(I)$ pour tout $k\leq 0$). De même, on dispose de son normalisé
$N(dR(D(I))) \in \fdg_{\WW_n(k)}$, construction qui définit un second $\s$-foncteur
$$N(dR(D(-))) : \mathcal{C}_n \longrightarrow \fdg_{\WW_n(k)}.$$
Or le morphisme naturel $A \to D(I)$ induit un morphisme sur les complexes de de Rham, 
qui est clairement compatible aux filtrations, et donc un morphisme de complexes filtrés
sur les normalisés correspondants
$$N(dR(A)) \longrightarrow N(dR(D(I))).$$
Ce morphisme est clairement fonctoriel en $A$ et donc définit une transformation naturelle
de $N(dR(-))$ vers $N(dR(D(-)))$.

D'après \cite[V-2.3.2]{bert}, tous les morphismes 
de transitions dans la direction simpliciale
$dR(D(I_0)) \longrightarrow dR(D(I_m))$
sont des quasi-isomorhismes filtrés, et $dR(D(I_0))$ est naturellement équivalent au complexe
de cohomologie cristalline de $\pi_0(A)$. Ainsi, la projection naturelle de complexes filtrés
$$dR(D(I_0)) \longrightarrow N(dR(D(I)))$$
est en réalité un quasi-isomorphisme filtré. On obtient ainsi un diagramme de complexes filtrés, 
fonctoriels en $A$
$$\xymatrix{N(dR(A)) \ar[r] & N(dR(D(I))) & \ar[l] dR(D(I_0))\simeq C^*_{cris}(\Spec\, A/S_n),}$$
et le second morphisme est une équivalence dans $\fdg_{\WW_n(k)}$. En passant aux complétés 
le long des filtrations on trouve un morphisme
$$C^*_{dR}(\Spec\, A/S_n)=\widehat{N(dR(A))} \longrightarrow 
\widehat{C}^*_{cris}(\Spec\, \pi_0(A)/S_n),$$
où le membre de droite est le complété du complexe de cohomologie cristalline de $\Spec\, \pi_0(A)$
pour la filtration de Hodge $PD$-adique. Ce morphisme est bien entendu fonctoriel 
en $A$. Cependant, il nous reste à voir que lorsque $\pi_0(A)$ est une $k$-algèbre lisse alors
le morphisme de complétion
$$C^*_{cris}(\Spec\, \pi_0(A)/S_n) \longrightarrow \widehat{C}^*_{cris}(\Spec\, \pi_0(A)/S_n)$$
est une équivalence de complexes de $\WW_n(k)$-modules. 

\begin{lem}\label{lapp}
Soit $A$ une $k$-algèbre lisse et $X=\Spec\, A$. Alors, le morphisme obtenu par complétion
le long de la filtration de Hodge PD-adique
$$C^*_{cris}(X/S_n) \longrightarrow \widehat{C}_{cris}^*(X/S_n)$$
est un quasi-isomophisme.
\end{lem}

\textit{Preuve du lemme.} On applique \cite[V-2.3.2]{bert} qui nous permet de calculer
la cohomologie cristalline en choisissant n'importe quel plongement $X \hookrightarrow Y$ avec $Y$ lisse
sur $S_n$. On prend pour $Y=\Spec\, B$ un relèvement lisse de $X$ à $S_n$ (qui existe
car $X$ est affine), et $X \hookrightarrow Y$
l'inclusion de la fibre spéciale. Dans ce cas, $A\simeq B/\pi$, avec $\pi$ une uniformisante dans
$\WW(k)$, et $B$ est sa propre enveloppe à puissance divisé de $(\pi)$. Ainsi, 
la cohomologie cristalline de $X$ sur $S_n$ se calcule simplement par le complexe de de Rham
de $Y$ au-dessus de $S_n$. La filtration de Hodge PD-adique sur ce dernier est finie (on utilise ici 
que $p\neq 2$, voir \cite[I-3.2.4]{bert}), et 
donc complète. \hfill $\Box$ \\

Nous avons donc construit le morphisme cherché: pour $A \in \mathcal{C}_n$, tel que
$\pi_0(A)$ soit une $k$-algèbre lisse et de plus $\pi_i(A)\simeq 0$ pour $i>0$, on trouve
le morphisme cherché
$$\psi : C^*_{dR}(\Spec\, \pi_0(A)/S_n) \longrightarrow 
\widehat{C}^*_{cris}(\Spec\, \pi_0(A)/S_n)\simeq C^*_{cris}(\Spec\, \pi_0(A)/S_n).$$

\end{appendix}

\bibliographystyle{alpha}
\bibliography{Chern-Fol-Bib}
\

\end{document}